\documentclass[12pt,reqno]{amsart}
\usepackage[cp1251]{inputenc}
\usepackage[english]{babel}
\usepackage{amssymb}
\usepackage{amsxtra}
\usepackage{mathrsfs}
\usepackage[all]{xy}
\usepackage{cite}
\usepackage{paralist}
\usepackage[hypertex]{hyperref}
\usepackage[active]{srcltx} 
\theoremstyle{plain}
\newtheorem{theorem}{Theorem}[section]
\newtheorem{lemma}[theorem]{Lemma}
\newtheorem{proposition}[theorem]{Proposition}
\newtheorem{corollary}[theorem]{Corollary}
\theoremstyle{definition}
\newtheorem{definition}[theorem]{Definition}
\newtheorem{remark}[theorem]{Remark}
\newtheorem{example}[theorem]{Example}
\newtheorem{question}[theorem]{Question}

\numberwithin{equation}{section}

\begin{document}
\title{Functions of class $C^\infty$ in non-commuting variables in the context of triangular Lie algebras}

\author{O.\,Yu.~Aristov}
\email{aristovoyu@inbox.ru}
\thanks{This work was supported by RFBR (grant no. 19-01-00447).}

\begin{abstract}
We construct a certain completion $C^\infty_\mathfrak{g}$ of the universal enveloping algebra of a triangular real Lie algebra $\mathfrak{g}$. It is a Fr\'echet-Arens-Michael algebra that consists of elements of polynomial growth and satisfies to the following universal property: every Lie algebra homomorphism from $\mathfrak{g}$ to a real Banach algebra all of whose  elements are of polynomial growth has an extension to a continuous homomorphism with domain~$C^\infty_\mathfrak{g}$. Elements of this algebra can be called functions of class $C^\infty$  in non-commuting variables. The proof is based on representation theory and employs an ordered $C^\infty$-functional calculus.
Beyond the general case, we analyze two simple examples. As an auxiliary material,  the basics of the general theory of algebras of polynomial growth are developed. We also consider local variants of the completion and obtain  a sheaf of non-commutative functions on the Gelfand spectrum of $C^\infty_\mathfrak{g}$ in the case when $\mathfrak{g}$ is nilpotent. In addition, we discuss the theory of holomorphic functions in non-commuting variables introduced by Dosi and
apply our methods to prove theorems strengthening some his results.
\end{abstract}

\maketitle

\markright{Non-commutative $C^\infty$-functions}

\begin{flushright}
\textit{For the memory of my parents,}
\\
\textit{Yurii and Liudmila Aristov}
\end{flushright}

\tableofcontents

\section{Introduction}
\label{s1}

Our purpose is to define an algebra of `infinitely smooth functions' on the non-commuta\-tive space corresponding to the universal enveloping algebra of a soluble real Lie algebra. We justify the definition by proving a universal property similar to that holding in the commutative case, namely, a theorem on $C^\infty$-functional calculus. Here we consider a calculus in a global sense, that is, the joint spectrum of the generators of the Lie algebra is not taken into account. (However, some results of the local theory, in particular, on sheaves of non-commutative algebras, are also included in this paper; see \S\,\ref{s5}.)

The method employed here leads to the goal when the Lie algebra is triangular and, as a simple example shows, it  is not applicable otherwise. We start from the following requirements for the desired algebra $C^\infty_\mathfrak{g}$ of `non-commutative infinitely smooth functions' associated with a triangular Lie algebra~$\mathfrak{g}$ over $\mathbb{R}$.

1.~In the abelian case, that is, when $\mathfrak{g}\cong \mathbb{R}^m$ for some $m\in\mathbb{N}$, the algebra $C^\infty_\mathfrak{g}$ must coincide with $C^\infty(\mathbb{R}^m)$, the algebra of real-valued infinitely differentiable functions on~$\mathbb{R}^m$.

2.~$C^\infty_\mathfrak{g}$ must be a completion of the universal enveloping algebra~$U(\mathfrak{g})$ and the homomorphism $U(\mathfrak{g})\to C^\infty_\mathfrak{g}$ must be injective if possible.

3.~The correspondence $\mathfrak{g}\mapsto C^\infty_\mathfrak{g}$ must be extended to a functor from  the category  of real Lie algebras (or, at least, the subcategory of triangular finite-dimensional ones) to a certain category of associative topological $\mathbb{R}$-algebras\footnote{We always assume that an associative algebra has an identity and a homomorphism of such algebras preserves the identities.} (of course, it must contain $C^\infty(\mathbb{R}^m)$ for every $m\in\mathbb{N}$).

For technical reasons, we are forced to abandon  the use of involutions (in algebras over $\mathbb{C}$) and work directly with algebras over reals\footnote{At the dawn of the theory of Banach algebras, the case of the field~$\mathbb{R}$ was considered on an equal footing with the case of~$\mathbb{C}$ and was more popular than nowadays. The fundamental results are the same in both cases; see \cite{1}.}.

In non-commutative differential geometry, certain  dense self-adjoint subalgebras of $C^*$-algebras are traditionally treated as algebras of `smooth functions'  (by analogy with $C^\infty(M)$ embedded in $C(M)$, where $M$ is a compact manifold). The algebras discussed here are not necessarily of this type. Moreover, we include  algebras with non-zero Jacobson radical in our framework. Assuming that the radical may be non-trivial, we expand the class of topological algebras under consideration.  This is the first of the two main ideas for the definition of `functions of class $C^\infty$ in non-commuting variables' proposed in this paper. Such an approach can be fruitful,  as shown by Dosi \cite{2}, who constructed a topological algebra of `formally radical entire functions' corresponding to a positively graded nilpotent complex Lie algebra. Some assertions generalizing Dosi's results are given in the last section, \S\,\ref{s6}.

To choose the category of topological algebras (see Requirement~3) it is necessary to select as constraints those features of $C^\infty(\mathbb{R}^m)$ that are essential for the problem under consideration. We use the following properties.

A.~It is an Arens--Michael algebra (that is, the topology can be determined by a family of submultiplicative seminorms).

B.~It is a Fr\'echet algebra (that is, it admits a countable family of seminorms determining the topology).

C.~For every $f\in C^\infty(\mathbb{R}^m)$ and every continuous submultiplicative seminorm $\|\,{\cdot}\,\|$ the function  $s\mapsto \|e^{isf}\|$, where $s\in\mathbb{R}$, is of polynomial growth (here $\|\,{\cdot}\,\|$ is extended to the complexification).

D.~For every $m$-tuple of elements of polynomial growth in a Banach algebra there exists a $C^\infty(\mathbb{R}^m)$-functional calculus.

Properties~A,~C and~D are essential in our considerations but Property~B is less important.

Now we consider functional calculi and elements of polynomial growth (in the sense of Property~C; see Definition~\ref{d2.1}) in more detail since their use is our second main idea.
Functional calculi of two types are traditionally employed in functional analysis: when the calculus is a homomorphism and when it is a linear mapping, not necessarily multiplicative (the latter is often used in the theory of non-commuting tuples of operators). Since we need both types, the term `multiplicative functional calculus' is used  in the first case and simply `functional calculus' in the second.  A so-called `ordered functional calculus' is a special case of the latter. We first consider it in a~purely algebraic form.

Let $\mathbb{R}[\lambda_1,\dots,\lambda_m]$ denote the algebra of polynomials in $m$ variables with real coefficients. If $B$ is an associative $\mathbb{R}$-algebra (not necessarily commutative) and $b_1,\dots,b_m\in B$, then the map taking a (commutative) monomial $\lambda_1^{\beta_1}\cdots \lambda_m^{\beta_m}$ to the (non-commutative) monomial $b_1^{\beta_1}\cdots b_m^{\beta_m}$  evidently extendes to a linear map
\begin{equation}
\label{eq1.1}
\Phi\colon \mathbb{R}[\lambda_1,\dots,\lambda_m]\to B.
\end{equation}
(Here we break the tradition of placing the factors in reverse order; see, for example, \cite[Appendix~1]{3}, but this is not essential.) This `non-commutative polynomial calculus' is the simplest example of an ordered calculus.

We are interested in conditions under which $\Phi$ can be extended from $\mathbb{R}[\lambda_1,\ldots,\lambda_m]$ to $C^\infty(\mathbb{R}^m)$. Denote the embedding of $\mathbb{R}[\lambda_1,\ldots,\lambda_m]$ in $C^\infty(\mathbb{R}^m)$ by $\iota$. Now let $B$ be a Banach or Arens-Michael algebra over~$\mathbb{R}$. We say that a continuous linear map $\theta\!:C^\infty(\mathbb{R}^m)\to B$ preserving unit is an \emph{ordered $C^\infty$-functional calculus} for $b_1,\ldots,b_m$ if the diagram
\begin{equation}
\label{eq1.2}
\begin{gathered}
\xymatrix@C=35pt{
\mathbb{R}[\lambda_1,\dots,\lambda_m]\ar[d]_{\iota}\ar[dr]^\Phi &
\\
C^\infty(\mathbb{R}^m)\ar[r]_{\theta} & B
}
\end{gathered}
\end{equation}
commutes. If an ordered $C^\infty$-functional calculus does exist, then it is unique (because  $\iota$ has dense range). If, in addition, $b_1, \ldots, b_m$ commute pairwise, then $\theta$ is obviously a homomorphism and then we call it a \emph{multiplicative $C^\infty$-functional calculus}. We use a similar terminology for other algebras, such as $C^n[c,d]$ and the non-commutative algebra $C^\infty_\mathfrak{g}$ defined below.

It is well known that, for a given element~$b$ of polynomial growth in a Banach algebra, one can use Fourier analysis to assign to each function $F\in C^\infty(\mathbb{R})$  an element $F(b)$  and thus obtain a $C^\infty$-functional calculus. The corresponding construction can be easily applied to the case of several pairwise commuting elements (see Theorem~\ref{t3.2}). In the case of $m$-dimensional Lie algebra~$\mathfrak{g}$,  we identify  $U(\mathfrak{g})$ with $\mathbb{R}[\lambda_1,\dots,\lambda_m]$ and aim to construct an ordered $C^\infty$-functional calculus similar to that in~\eqref{eq1.2} and, moreover, multiplicative, that is, a homomorphism. (Results on $C^\infty$-functional calculi in the form necessary for what follows are contained in~\S\,\ref{s3}. Proofs are summarized there for completeness since the author was unable to find them in the available literature. Certain necessary modifications of the results for the case when some of the elements are nilpotent are also included in~\S\,\ref{s3}.)

Since the property of being an element of polynomial growth is not preserved, in generally, neither under passage to a limit nor under taking algebraic combinations of non-commuting elements (see Remark~\ref{r2.7}), we consider Banach or Arens--Michael $\mathbb{R}$-algebras  \textbf{all of whose elements} have polynomial growth and call them \textit{algebras of polynomial growth}. The class of such algebras is basic for us and their properties are discussed in~\S\,\ref{s2}. The main example is the algebra of $n$ times continuously differentiable functions with values in triangular real matrices. (Here it is essential that we work with the field~$\mathbb{R}$.)

Linear operator of polynomial growth have been known for a long time and have been studied in  the theory of generalized scalar operators initiated in the paper~\cite{4} of  Foia\c{s} in~1960. However, as far as the author knows, non-commutative Banach algebras \textbf{entirely} consisting of elements satisfying this condition have not been previously considered.

The image of a finite-dimensional Lie $\mathbb{R}$-algebra under a homomorphism to a projective limit of Banach  $\mathbb{R}$-algebras of polynomial growth is a triangular finite-dimensional Lie algebra (Corollary~\ref{c4.2}). Thus, if we want $U(\mathfrak{g})\to C^\infty_\mathfrak{g}$ to be injective (see Requirement~2), then we have to assume that~$\mathfrak{g }$ is triangular. Note that in this case an ordered $C^\infty$-functional calculus can be treated as an analytic version of the theorem on the existence of a PBW-basis in $U(\mathfrak{g})$.  A more general situation, where $\mathfrak{g}$ is finitely generated but not necessarily finite dimensional, is considered in~\cite{5}, which is a follow-up to this article.

The main results of the paper are stated in \S\,\ref{ss4.2}. Specifically, suppose that $\mathfrak{g}$ is triangular, $e_{k+1},\dots,e_m$ is a linear basis of its commutant and $e_1,\dots, e_k$ completes it to a linear basis of~$\mathfrak{g}$. The first main result, Theorem~\ref{t4.3}, asserts that the multiplication in $U(\mathfrak{g})$ can be extended by continuity to the projective tensor product
$$
C^\infty_\mathfrak{g}:=C^\infty(\mathbb{R}^k)\mathbin{\widehat{\otimes}} \mathbb{R}[[e_{k+1},\dots,e_m]]
$$
of the function space and the formal power series space. The algebra defined in this way has the required properties. The proof makes essential use of both representation theory and the structural theory of soluble Lie algebras. Note that even the continuity of the multiplication in $C^\infty_\mathfrak{g}$ is far from obvious (sf. a `holomorphic' version in the nilpotent case in~\cite{2}). Our argument includes many technical details but the main idea can be clearly outlined already for the two simplest examples, the two-dimensional non-abelian algebra and the three-dimensional Heisenberg algebra. They are discussed separately in~\S\,\ref{ss4.4}.

In the second main result, Theorem~\ref{t4.4}, we prove a universal property for $C^\infty_\mathfrak{g}$ (in other words, the existence of a multiplicative $C^\infty$-functional calculus mentioned above). It follows, in particular, that $C^\infty_\mathfrak{g}$ is independent of the choice of a basis. Note that Theorem~\ref{t4.4} can be interpreted as a statement that $C^\infty_\mathfrak{g}$ is the envelope of the associative algebra $U(\mathfrak{g})$ with respect to the class of Banach algebras of polynomial growth (sf. the Arens-Micheal envelope, which uses the class of \textbf{all} Banach algebras \cite{6}). This envelope is considered in more detail in~\cite{5}.

In the last sections of the article, we discuss two topics closely related to its main subject and obtain a number of results based on auxiliary statements in \S\,\ref{ss4.5} and their simple modifications. The local version of the theory, with definition of algebras of non-commutative functions but without considering the functional calculus (which is expected in the future), is examined in~\S\,\ref{s5}. In particular, we construct sheaves of non-commutative functions  in the case of a nilpotent Lie algebra.

The algebra $C^\infty_\mathfrak{g}$ defined in this paper and its local variants are analogues of the algebra of `formally radical holomorphic functions' studied by Dosi in~\cite{2} and \cite{7}. The technique employed here  makes it possible to obtain assertions similar to that in the $C^\infty$-theory and, in particular, to improve some of Dosi's results; see~\S\,\ref{s6}.

In the future, in a separate article, the author hopes to study the relationship between the algebras considered in this paper and the algebras of `smooth functions' that arise in the theory of $C^*$-algebras and non-commutative geometry in the spirit of Connes (for variants of the theorem on functional calculus see \cite[Proposition~6.4]{8}, \cite[Proposition 22]{9} and \cite[Proposition 2.8]{10}).

\section{Elements and algebras of polynomial growth}
\label{s2}

\subsection{Elements of polynomial growth}

\begin{definition}
\label{d2.1}
An element $b$ of a complex Arens-Michael (in, particular, Banach) algebra~$B$  is of
\textit{polynomial growth}\footnote{This definition treats elements with real spectrum. An alternative definition for elements with spectrum contained in the unit circle includes a similar estimate for $\|b^n\|$, where $n\in\mathbb{Z}$, but we do not need it.} if for every continuous submultiplicative seminorm $\|\,{\cdot}\,\|$ on~$B$ there are $K>0$ and $\alpha\ge0$ such that
\begin{equation}
\label{eq2.1}
\|e^{isb}\|\le K (1+|s|)^{\alpha} \quad \text{for all }s\in \mathbb{R}.
\end{equation}
An element of a real algebra Arens-Michael is of \textit{polynomial growth}
if it is of polynomial growth in the complexification.
\end{definition}

Note that the construction of the complexification of a real Banach algebra (see \cite[Chapter~I, \S\,3, p.\,5--9]{1} and~\cite[\S\,2.1]{11}) can be easily carried over to the case of an Arens--Michael algebra and $e^{isb}$ is well defined since it is always possible to substitute an element of an Arens--Michael algebra into a given entire function (in other words,~$b$ admits a holomorphic functional calculus on~$\mathbb{C}$).

Our terminology is close to that used in~\cite{12}. Alternative names are `element of slow growth' (`\/\`{e} croissance lente')~\cite{13}, `generating element' \cite{14} and `generalized scalar element' \cite{15}, \cite{16}. Note that the last term is usually used in operator theory in a more general sense.

Every element of polynomial growth has real spectrum\footnote{By the spectrum of an element of an Arens--Michael algebra we always mean the spectrum in the complexification.}, that is, containing in~$\mathbb{R}$ (see Proposition~\ref{p2.2}). The converse does not hold in general but it is true in the particular case of interest to us, namely, for the algebra of infinitely smooth functions with values in triangular matrices. To verify this we use the following assertion.

\begin{proposition}
\label{p2.2}
Let $b$ be an element of a Banach algebra~$B$ over $\mathbb{C}$. The following conditions are equivalent.

\emph{(1)}~$b$ is of polynomial growth.

\emph{(2)}~There are a non-trivial interval $[c,d]$ in~$\mathbb{R}$ and $n\in\mathbb{N}$ such that there exists a multiplicative  functional calculus $C^n[c, d]\to B$ for~$b$.

\emph{(3)}~The spectrum of $b$ is contained in~$\mathbb{R}$ and there are $C,\gamma>0$ such that
$$
\|(b-\lambda)^{-1}\|\le C(1+|{\operatorname{Im} \lambda}|)^{-\gamma}
$$
when~$\operatorname{Im}\lambda\ne0$.
\end{proposition}

It is easy to see that a functional calculus for $\mathbb{R}$-valued functions of class $C^n$ extends to $\mathbb{C}$-valued functions. For the latter case, the proof repeats verbatim the argument for the Colojoar\v{a}--Foia\c{s} theorem (as given in~\cite[Theorem 1.5.19]{17}), which uses only the structure of a Banach algebra on the set of operators; see also \cite[Theorem~1]{13}.

We need primary properties of elements of polynomial growth.

\begin{proposition}
\label{p2.3}
Let $\varphi\colon A\to B$ be a continuous homomorphism of Arens--Michael algebras. If $a\in A$ is of polynomial growth, then so is~$\varphi(a)$.
\end{proposition}

\begin{proof}
Let $\|\,{\cdot}\,\|$ be a continuous submultiplicative seminorm on~$B$. Since $\varphi$ is continuous, there are a continuous submultiplicative seminorm $\|\,{\cdot}\,\|'$ on~$A$ and $C>0$ such that $\|\varphi(a)\|\le C\|a\|'$ for every $a\in A$. If $a$ is of polynomial growth, then, by definition, there are $K>0$ and $\alpha\ge0$ such that $\|e^{isa}\|'\le K (1+|s|)^{\alpha} $ for all $s\in\mathbb{R}$. Since $e^{is\varphi(a)}=\varphi(e^{isa})$, it follows that $\varphi(a)$ is also of polynomial growth.
\end{proof}

Recall that an element $b$ of a Banach algebra is said to be \textit{topologically nilpotent} if $\lim_{n\to\infty}\|b^n\|^{1/n} =0$.

\begin{proposition}
\label{p2.4}
If an element of a Banach algebra is topologically nilpotent and of polynomial growth, then it is
nilpotent.
\end{proposition}

This fact is well known for operators in the Banach space; see, for example, \cite[Proposition 1.5.10]{17}. The proof for elements of a Banach algebra is identical.

\begin{proposition}
\label{p2.5}
In an Arens--Michael algebra, the spectrum of every element of polynomial growth  is contained in~$\mathbb{R}$.
\end{proposition}

\begin{proof}
The spectrum of every element of an Arens--Michael algebra is the union of the spectra of all corresponding elements in the concomitant Banach algebras \cite[Chapter~5, Russian p.\,280, Corollary~2.12]{6}. It remains to note that
by Proposition~\ref{p2.3}  these elements are of polynomial growth and so by the implication (1)$\Rightarrow$(3) in Proposition~\ref{p2.2} their spectra are contained in~$\mathbb{R}$.
\end{proof}

\subsection{Algebras of polynomial growth}
\label{ss2.2}

We now single out the main class of algebras used in this paper.

\begin{definition}
\label{d2.6}
An Arens-Michael $\mathbb{R}$-algebra is said to be of  \textit{polynomial growth} if all its elements are of polynomial growth.
\end{definition}

First we note that in the case of the field~$\mathbb{C}$ it is pointless to discuss algebras all of whose elements are of polynomial growth since even the imaginary unit~$i$ does not satisfy this condition. Therefore we consider only algebras over~$\mathbb{R}$.
\begin{remark}
\label{r2.7}
Considering algebras all of whose elements are of polynomial growth, we have to impose additional restrictions on topological generators. Indeed, it is easy to see that any algebraic combination of commuting elements of polynomial growth also has polynomial growth (cf. \cite[Corollary 1.5.20]{17} or \cite[Chapter~4, p.\,106, Corollary 3.4]{18}). However, this is not the case for non-commuting elements (counterexamples can be easily found among matrices of order~$2$ with real eigenvalues). Moreover, the property of having polynomial growth is not preserved in general under passage to the limit, even in the commutative case. To verify this, consider the Banach algebra $A(\mathbb{T})_\mathbb{R}$ of absolutely convergent real Fourier series. The functions $\sin$ and $\cos$ are of polynomial growth (this  follows, for example, from  \cite[Chapter~VI, \S\,2, \S\,3, Russian p.\,93--95]{19}). It is obvious that the set of their algebraic combinations is dense in $A(\mathbb{T})_\mathbb{R}$ and so $\sin$ and $\cos$ are topological generators of $A(\mathbb{ T})_\mathbb{R}$. However, the assumption that every element of $A(\mathbb{T})_\mathbb{R}$ is of polynomial growth contradicts Katznelson's theorem, which states that a function defined on an interval and acting on~$A(\mathbb{T})_\mathbb{R}$ is analytic (see \cite[Chapter~VI, \S\,6, Russian p.\,102]{19}). Indeed, if $A(\mathbb{T})_\mathbb{R}$ is of polynomial growth, then for each of its elements there  exists a $C^n$-calculus for some $n\in\mathbb{ N}$ (see Proposition~\ref{p2.2}), which means that any function of class $C^\infty$ acts on $A(\mathbb{T})_\mathbb{R}$. This clearly contradicts the fact that there are non-analytic $C^\infty$-functions.
\end{remark}

The following theorem is well known and the proof is given for completeness. We denote the Jacobson radical of $B$ by $\operatorname{Rad} B$.

\begin{theorem}
\label{t2.8}
Let $B$ be a Banach algebra over~$\mathbb{R}$ such that the spectrum of every element is contained in~$\mathbb{R}$ (in particular, an algebra all of whose elements are of polynomial growth). Then $B/\operatorname{Rad} B$ is commutative.
\end{theorem}

\begin{proof}
The spectrum of every element of $B/\operatorname{Rad} B$ is also contained in~$\mathbb{R}$ since it cannot be increased by taking a quotient. By a result of Kaplansky \cite[Theorem~4.8]{20},  a semisimple Banach algebra over~$\mathbb{R}$ is commutative when $1+x^2$ is invertible for every  element~$x$. Note that this condition is equivalent to the property that the spectrum of every element is contained in~$\mathbb{R}$. Since $B/\operatorname{Rad} B$ semisimple  \cite[Chapter~2, \S\,3, p.\,56]{1}, this means that $B/\operatorname{Rad} B$ is commutative.
\end{proof}

\begin{proposition}
\label{p2.9}
Let $B$ be a Banach $\mathbb{R}$-algebra of polynomial growth. Then $\operatorname{Rad} B$ is nilpotent \emph{(}that is, there is $n\in\mathbb{N}$ such that $r_1\cdots r_n=0$ for arbitrary $r_1,\dots ,r_n\in\operatorname{Rad} B$\emph{)}.
\end{proposition}

\begin{proof}
Every  element in the radical of a Banach algebra is topologically nilpotent. Therefore, by Proposition~\ref{p2.4}, every element of $\operatorname{Rad} B$ is nilpotent. As Grabiner  noted \cite{21}, it follows from the Dubnov--Ivanov theorem
(also known as the Nagata--Higman theorem) that, in the case of a Banach algebra, this condition implies that $\operatorname{Rad} B$ is nilpotent.
\end{proof}

\begin{remark}
\label{r2.10}
Note that every Banach $\mathbb{R}$-algebra of polynomial growth satisfies a polynomial identity. Indeed,  Proposition~\ref{p2.9} implies that there is $n$ such that $(\operatorname{Rad} B)^n=0$. By Theorem~\ref{t2.8}, we have $[a,b]\in \operatorname{Rad} B$ and so $[a,b]^n=0$ for every $a,b\in B$.
\end{remark}

\begin{proposition}
\label{p2.11}
{\rm(A)}~The class of Banach $\mathbb{R}$-algebras of polynomial growth is stable under passage to closed subalgebras and finite Cartesian products. The class of Arens--Michael $\mathbb{R}$-algebras of polynomial growth is stable under passage to  closed subalgebras and arbitrary Cartesian products, and, as a consequence, projective limits.

{\rm(B)}~An Arens--Michael $\mathbb{R}$-algebra is a projective limit of Banach algebras of polynomial growth if and only if it is isomorphic to a closed subalgebra of a product of Banach algebras of polynomial growth.
\end{proposition}

\begin{proof}
Part~(A) follows directly from the definitions.

The necessity in Part~(B) follows from the standard construction of the projective limit; see, for example, \cite[Chapter~III, \S\,2, p.\,84, Lemma~2.1]{22}. To prove the sufficiency note that a closed subalgebra~$A$ of a product of Banach algebras of polynomial growth has a base of neighbourhoods of~$0$ consisting of barrels stable under multiplication and the Banach algebra corresponding to each element of this base is of polynomial growth. By \cite[Chapter~III, \S\,3, p.\,88, Theorem~3.1]{22}, the algebra~$A$ is the projective limit of the Banach algebras corresponding to the elements of the base.
\end{proof}

Denote by $\mathrm{T}_p$ the algebra of upper triangular (including the diagonal) real matrices of order~$p$, where $p\in\mathbb{N}$. Our main model example of an algebra of polynomial growth is the algebra $C^\infty(V,\mathrm{T}_{p})$ of infinitely differentiable $\mathrm{T}_p$-valued functions on an open subset~$V$ of~$\mathbb{R}^m$. We fix a submultiplicative norm $\|\,{\cdot}\,\|_p$ on~$\mathrm{T}_p$ (for example, we can take the operator norm in $p$-dimensional Euclidean space). Then the topology on $C^\infty(V,\mathrm{T}_{p})$ is determined by the following family of seminorms:
\begin{equation}
\label{eq2.2}
\|f\|_{p,K,n}:=\sum_{\beta\in \mathbb{Z}_+^m,\,|\beta|=n}\|f^{(\beta)}\|_{p,K,0},\qquad
\text{where}\quad \|f\|_{p,K,0}:=\max_{x\in K}\|f(x)\|_p,
\end{equation}
$n\in\mathbb{Z}_+$, the partial derivative of vector-valued function $f\in C^\infty(V,\mathrm{T}_{p})$ is denoted by~$f^{(\beta)}$ and $K$ is a compact subset of~$V$. (Here we put $|\beta|:=\beta_1+\cdots+ \beta_m$ for $\beta=(\beta_1,\dots, \beta_m)\in\mathbb{Z}_+^m$.) It is easy to deduce from Leibniz's rule for the derivative of a composite of a bilinear map and vector-valued functions that the seminorms $\sum_{n=0}^q\|\,{\cdot}\,\|_{p,K,n}/n!$, $q\in\mathbb{Z}_+$, are submultiplicative. Since the family of all such seminorms is equivalent to the initial one, we see that $C^\infty(V,\mathrm{T}_{p})$ is an Arens--Michael algebra.

We also consider Banach algebras $C^n(K,\mathrm{T}_{p})$, where $n\in\mathbb{Z}_+$ and $K$ is a compact subset of~$\mathbb{R}^m$ with dense interior\footnote{Generally speaking, the requirement of having dense interior is redundant. We use it to avoid the subtleties associated with defining the class $C^n$ in the general case.}.

\begin{theorem}
\label{t2.12}
Let $p\in\mathbb{Z}_+$, $m\in \mathbb{N}$ and~$V$ be an open subset of~$\mathbb{R}^m$. Then
$C^\infty(V,\mathrm{T}_{p})$ is a projective limit of Banach algebras of polynomial growth.
\end{theorem}

Consider first the commutative case (that is, when $p=1$).

\begin{proposition}
\label{p2.13}
Let $K$ be a compact subset of~$\mathbb{R}^m$ with dense interior and $n\in\mathbb{Z}_+$. Then $C^n(K)$ is a Banach algebra of polynomial growth.
\end{proposition}

\begin{proof}
We present two proofs here: the first proof uses Proposition~\ref{p2.2} and the second is based directly on the definition.

\textsl{First proof}. Take $f\in C^n(K)$ and an interval $[c,d]$ containing the range of~$f$ and note that~$f$ induces a homomorphism $C^n[c,d]\to C^n(K)\colon g \mapsto (t \mapsto g(f(t)))$, which is a multiplicative functional calculus. It remains to apply Part~(2) of Proposition~\ref{p2.2}.

\textsl{Second proof}. Here and below, in the commutative case (that is, when $p=1$), we denote the seminorms $\|\,{\cdot}\,\|_{1,K,n}$ defined in~\eqref{eq2.2} by $|\,{\cdot}\,|_{K,n}$. Thus the topology on $C^n(K)$ is determined by  $|\,{\cdot}\,|_{K,n}$. It is easy to see that for every $f\in C^n(K)$ and $\beta\in\mathbb{Z}_+^m$ the partial derivative $(e^{isf})^{(\beta)} $ is equal to the product of $e^{isf}$ and a polynomial in~$s$ of degree $|\beta|$ with coefficients that are polynomials in partial derivatives  of~$f$. Since~$f$ takes only real values, we have $|e^{isf}|_{K,0}=1$. Hence $|e^{isf}|_{K,n}$ can be majorized by a polynomial in~$|s|$ of degree~$n$ and so $f$ is of polynomial growth.
\end{proof}

Now we  deduce Theorem~\ref{t2.12} from a general result. Recall first that an extension
$$
0\leftarrow A \leftarrow \mathfrak{A} \leftarrow I\leftarrow 0
$$
of a Banach algebra~$A$ by~$I$ is said to be \textit{nilpotent } if~$I$ is a closed nilpotent
ideal of~$\mathfrak{A}$ and \textit{splits} if there is a continuous homomorphism $A\to \mathfrak{A}$ that is right inverse for~$\mathfrak{A}\to A$. (For more details on extension of Banach algebras see \cite{6} or~\cite{23}.)

\begin{theorem}
\label{t2.14}
If an extension of a Banach $\mathbb{R}$-algebra of polynomial growth splits and is nilpotent, then it is also an algebra of polynomial growth.
\end{theorem}

It is well known that not every nilpotent extension splits even in the case when $A=C^n[0,1]$ and~$I$ is finite dimensional; see, for example,  \cite [\S\,5, pp.\,84--87]{23} for $\mathbb{C}$-algebras.

\begin{proof}[Proof of Theorem~\ref{t2.14}]
Let an extension $0\leftarrow A \leftarrow \mathfrak{A} \leftarrow I\leftarrow 0$ of Banach $\mathbb{R}$-algebras split and be nilpotent and let $A$~be an algebra of polynomial growth. Fix a splitting continuous homomorphism $\rho\colon A\to \mathfrak{A}$.

Let~$b\in \mathfrak{A}$. Then $b=d+t$, where $t\in I$ and $d=\rho(a)$ for some $a\in A$.
Hence $d$ is of polynomial growth by Proposition~\ref{p2.3}.

Suppose that $\lambda\in\mathbb{C}\setminus \mathbb{R}$. Since the spectrum of~$d$ in the complexification of $\mathfrak{A}_\mathbb{C}$ is contained in~$\mathbb{R}$, the element $d-\lambda$ is invertible. Then $(d-\lambda)^{-1}t$ is well defined and belongs to~$I$. Therefore it is nilpotent and so $1+ (d-\lambda)^{-1}t$ is invertible. From the obvious equality $d-\lambda+t=(d-\lambda)(1+ (d-\lambda)^{-1}t)$ we obtain that $b-\lambda$ is invertible and
$$
(b-\lambda)^{-1}=(1+ (d-\lambda)^{-1}t)^{-1}(d-\lambda)^{-1}.
$$
Since $I$ is nilpotent, there is $p\in\mathbb{N}$ such that $((d-\lambda)^{-1}t)^p=0$. It follows that
$$
(1+ (d-\lambda)^{-1}t)^{-1}= \sum_{j=0}^{p-1} ((\lambda-d)^{-1}t)^{j}.
$$
Hence,
$$
\|(b-\lambda)^{-1}\|\le\sum_{j=0}^{p-1} \|(d-\lambda)^{-1}\|^{j+1} \|t\|^j.
$$

Applying the implication $(1)\Rightarrow(3)$ in Proposition~\ref{p2.2} to $d$, we conclude that there exist $C,\gamma>0$ such that for every $\lambda$ satisfying $\operatorname{Im} \lambda\ne0$ the inequality
$$
\|(d-\lambda)^{-1}\|\le C(1+|{\operatorname{Im} \lambda}|)^{-\gamma}
$$
holds.

Thus we have for $\operatorname{Im} \lambda\ne0$ that
$$
\|(b-\lambda)^{-1}\|\le C'(1+|{\operatorname{Im} \lambda}|)^{-\gamma'}
$$
with $C'$ and $\gamma'$  depending only on~$p$ and~$\|t\|$. It follows from the converse implication $(3)\Rightarrow(1)$ in Proposition~\ref{p2.2} that~$b$ is of  polynomial growth in $\mathfrak{A}_\mathbb{C}$ and then by the definition also in~$\mathfrak{A}$.
\end{proof}

In what follows we use the sign $\mathbin{\widehat{\otimes}}$ for the complete projective product of locally convex spaces. Recall that if $B_1$ and $B_2$ are (Arens--Michael) Banach algebras, then so is $B_1\mathbin{\widehat{\otimes}} B_2$.

\begin{proof}[Proof of Theorem~\ref{t2.12}]
The algebra $C^\infty(V,\mathrm{T}_{p})$ is the projective limit of algebras of the form $C^n(K,\mathrm{T}_{p})$, where~$K$ is a compact subset of~$\mathbb{R}^m$ and $n\in\mathbb{N}$. Moreover, we can assume that $K$ has dense interior. Therefore, it suffices to show that every algebra of this type is of polynomial growth. Note that $C^n(K,\mathrm{T}_{p})$ is topologically isomorphic to the projective tensor product  $C^n(K)\mathbin{\widehat{\otimes}} \mathrm{T}_{p}$ (because $\mathrm{T}_{p}$ is finite dimensional).

Denote by $I_0$ the ideal of~$\mathrm{T}_{p}$ consisting of all the matrices with zero diagonal entries. Then the extension
$$
0\leftarrow C^n(K)^p \leftarrow C^n(K)\mathbin{\widehat{\otimes}} \mathrm{T}_{p} \leftarrow C^n(K)\mathbin{\widehat{\otimes}} I_0\leftarrow 0
$$
splits and is nilpotent
since so is the extension
$$
0\leftarrow \mathbb{R}^p \leftarrow \mathrm{T}_{p} \leftarrow I_0\leftarrow 0.
$$
Thus $C^n(K)^p$ is an algebra of polynomial growth by Propositions~\ref{p2.11} and~\ref{p2.13}. Finally we apply Theorem~\ref{t2.14}.
\end{proof}


\section{Ordered calculus}
\label{s3}

The following theorem is not new (as well as its consequence, Theorem~\ref{t3.2}). A one-dimensional variant~\cite{13} and a multi-dimensional variant for a certain class of smooth symbols \cite[p.\,271]{3} are well known. The result for the class $C^\infty$ is stated in~\cite[Russian p.\,888]{14} but with details omitted. We give a summary of a proof. (Here we denote by $\mathcal{S}(\mathbb{R}^m)$ the Schwartz space of rapidly decreasing infinitely differentiable functions.)

\begin{theorem}
\label{t3.1}
Let  $\mathbf{b}=(b_1, \dots, b_m)$ be an ordered tuple of elements of polynomial growth of a Banach algebra~$B$. For $F\in C^\infty(\mathbb{R}^m)$ take $f\in \mathcal{S}(\mathbb{R}^m)$ coinciding with~$F$ on the parallelepiped
\begin{equation}
\label{eq3.1}
P:=[-\|b_1\|,\|b_1\|]\times\dots\times [-\|b_m\|,\|b_m\|],
\end{equation}
and put
\begin{equation}
\label{eq3.2}
f(\mathbf{b}):=\frac1{(2\pi)^m}\int_{\mathbb{R}^m} \widehat f(\mathbf{s}) \exp(i s_1b_1)\cdots \exp(i s_mb_m)\,d\mathbf{s},
\end{equation}
where $\widehat f$ is the Fourier transform of~$f$. Then $\theta(F):=f(\mathbf{b})$ is independent of the choice of $f$ and $\theta$ is an ordered $C^\infty$-functional calculus for~$\mathbf{b}$.
\end{theorem}

(For  ordered $C^\infty$-functional calculi see \S\,\ref{s1}.)

\begin{proof}
First, note that the integral in~\eqref{eq3.2} is well defined for every $f\in \mathcal{S}(\mathbb{R}^m)$. This assertion can be proved in the same way as in the one-dimensional case; see, for example, \cite[Lemma 1.5.18]{17}. Thus $f\mapsto f(\mathbf{b})$ is a $B$-valued distribution of class $\mathcal{S}'$.

Second, we claim that the support of this distribution is compact and contained in~$P$. Here we use the same argument  as for the Weyl calculus (see, for example, \cite[lemma~8.4]{24}).
Put
$$
E(\mathbf{z}):=(2\pi)^{-m}\exp(i z_1b_1)\cdots \exp(i z_mb_m),\qquad \mathbf{z}=(z_1, \dots, z_m)\in\mathbb{C}^m.
$$
It is obvious that $E$ is an entire $B$-valued function in  $m$ variables. Write $\mathbf{z}=\mathbf{x}+i\mathbf{y}$, where $\mathbf{x},\mathbf{y}\in \mathbb{R}^m$.

Since $b_1,\dots,b_m$ are of polynomial growth, there are $K>0$ and $\alpha\ge0$ such that
$$
\|E(\mathbf{z})\|\le K(1+ |\mathbf{x}|)^{\alpha} \exp\biggl(\sum_j \|b_j\|\, |y_j|\biggr)
$$
for every $\mathbf{z}\,{\in}\,\mathbb{C}^m$.
Let~$\eta$ be an arbitrary continuous linear functional on~$B$. By the Paley--Wiener--Schwartz theorem (as it stated in~\cite[Theorem~7.3.1]{25}),  the function $\mathbf{z}\mapsto \eta E(\mathbf{z})$ is the Fourier transform of a distribution~$u$ with support contained in~$P$. In particular, if $f\in \mathcal{S}(\mathbb{R}^m)$, then $\langle \eta E, \widehat f\rangle = \langle u,f\rangle $ and, moreover, \eqref{eq3.2} implies that $\eta f(\mathbf{b})=\langle \eta E, \widehat f\rangle $. Thus  $f\mapsto f(\mathbf{b})$ is a  $B$-valued distribution with compact support.

The support of the distribution is contained in~$P$. So if $F\in C^\infty(\mathbb{R}^m)$ and $f$ is a function in $\mathcal{S}(\mathbb{R}^m) $ such that $f=F$ on~$P$, then $f(\mathbf{b})$ is independent of the choice of~$f$.
Thus the formula $\theta(F)\,{:=}\,f(\mathbf{b})$ 
defines a continuous linear map $\theta\colon C^\infty(\mathbb{R}^m)\,{\to}\, B$
independent of the choice of~$f$.

Finally, we verify that $\theta$ is an ordered $C^\infty$-functional calculus. It suffices to show that  $\theta(x_1^{\alpha_1}\cdots x_m^{\alpha_m})=b_1^{\alpha_1}\cdots b_m^{\alpha_m}$ for every multi-index~$\alpha$. We can replace the function $\mathbf{x}\mapsto x_1^{\alpha_1}\cdots x_m^{\alpha_m}$ by $f\in \mathcal{S}(\mathbb{R}^m)$ that coincides with it on~$P$. Specifically, take for every~$j$ a function $\varphi_j$ of class $C^\infty$ on $\mathbb{R}$  that equals~$1$ on~$[-\|b_j\|,\|b_j\|]$ and have compact support and then put
$$
f(\mathbf{x}):=x_1^{\alpha_1}\cdots x_m^{\alpha_m}\varphi_1(x_1)\cdots \varphi_m(x_m).
$$
Since $f\in\mathcal{S}(\mathbb{R}^m)$, we have $\theta(x_1^{\alpha_1}\cdots x_m^{\alpha_m})=f(\mathbf{b})$. Then the integral in~\eqref{eq3.2} is a product of one-dimensional integrals and thus we can reduce the assertion to the case of one variable, for which it is known. Indeed, the functional $C^n$-calculus in Part~(2) of Proposition~\ref{p2.2} can be defined by the one-dimensional version of~\eqref{eq3.2} (see \cite[Theorem~1, Lemma~1]{13} or \cite[Lemmas 1.5.16, 1.5.18]{17}).
\end{proof}

In the particular case when $b_1, \dots, b_m$ commute pairwise, the exponentials also commute and we immediately obtain from Theorem~\ref{t3.1} the following result on a multiplicative $C^\infty$-calculus.

\begin{theorem}
\label{t3.2}
Let $\mathbf{b}=(b_1, \dots, b_m)$ be a tuple of pairwise commuting elements of polynomial growth in a Banach algebra~$B$. For $F\in C^\infty(\mathbb{R}^m)$ take $f\in \mathcal{S}(\mathbb{R}^m)$ coinciding with~$F$ on the parallelepiped  \footnote{In fact, the support of $\theta$ coincides with the joint spectrum defined in almost any reasonable sense \cite{26}, but we do not need this fact.} defined in~\eqref{eq3.1} and put
\begin{equation}
\label{eq3.3}
f(\mathbf{b}):=\frac1{(2\pi)^m}\int_{\mathbb{R}^m} \widehat f(\mathbf{s}) \exp(i s_1b_1+\cdots +i s_mb_m)\,d\mathbf{s}.
\end{equation}
Then $\theta(F):=f(\mathbf{b})$ is independent of the choice of~$f$ and determines a multiplicative functional calculus, that is, a continuous homomorphism
$$
\theta\colon C^\infty(\mathbb{R}^m) \to B
$$
such that $\theta(t_j)=b_j$ for all $j=1,\dots,m$.
\end{theorem}

In what follows, we need a strengthening of Theorem~\ref{t3.1} in the case when some of $b_1, \dots, b_m$ are nilpotent. (Note that every nilpotent element is of polynomial growth.)

Denote by $\mathbb{R}[[\lambda_{k+1},\dots,\lambda_m]]$ the linear space of formal power series in variables $\lambda_{k+1},\dots,\lambda_m$ with real coefficients. It can be treated as the quotient algebra of $C^\infty(\mathbb{R}^{m-k})$ by the closed ideal~$I$ consisting of functions all of whose  partial derivatives vanishing at the origin (this follows from Borel's theorem; see, for example, \cite[Chapter~I, p.\,18, Theorem 1.3]{27}). The corresponding topology is the product topology on the product of a countable number of copies of~$\mathbb{R}$. Since $C^\infty(\mathbb{R}^k)$ is a nuclear Fr\'echet space, the linear map
$$
C^\infty(\mathbb{R}^k)\mathbin{\widehat{\otimes}} I\to C^\infty(\mathbb{R}^k)\mathbin{\widehat{\otimes}} C^\infty(\mathbb{R}^{m-k})
$$
is topologically injective and
\begin{equation}
\label{eq3.4}
C^\infty(\mathbb{R}^k)\mathbin{\widehat{\otimes}} \mathbb{R}[[\lambda_{k+1},\dots,\lambda_m]]\cong (C^\infty(\mathbb{R}^k)\mathbin{\widehat{\otimes}} C^\infty(\mathbb{R}^{m-k}))/(C^\infty(\mathbb{R}^k)\mathbin{\widehat{\otimes}} I);
\end{equation}
see, for example, \cite[Theorem~A1.6]{28}.

\begin{theorem}
\label{t3.3}
Let $\mathbf{b}=(b_1, \dots, b_m)$ be an ordered tuple of elements of an Arens-Michael  \emph{(}in particular, Banach\emph{)} algebra $\mathbb{R}$-algebra~$B$. Suppose that $b_1,\dots,b_k$ \emph{(}$k\le m$\emph{)} are of polynomial growth and $b_{k+1},\dots,b_m$ are nilpotent. For $F\in C^\infty(\mathbb{R}^k)\mathbin{\widehat{\otimes}} \mathbb{R}[[\lambda_{k+1},\dots,\lambda_m]]$ and some its pre-image $f\in \mathcal{S}(\mathbb{R}^m)$ we define $f(\mathbf{b})$ as in~\eqref{eq3.2}. Then $\sigma(F):=f(\mathbf{b})$ is independent of the choice of~$f$ and determines a linear map
$$
\sigma\colon C^\infty(\mathbb{R}^k)\mathbin{\widehat{\otimes}} \mathbb{R}[[\lambda_{k+1},\dots,\lambda_m]]\to B
$$
extending the map in~\eqref{eq1.1}.
\end{theorem}

By the Weierstrass approximation theorem for continuously differentiable functions, the space $\mathbb{R}[\lambda_1,\dots,\lambda_m]$ is dense in $C^\infty(\mathbb{R}^k)$. Then  the range of the map
$$
\mathbb{R}[\lambda_1,\dots,\lambda_m]\to C^\infty(\mathbb{R}^k)\mathbin{\widehat{\otimes}} \mathbb{R}[[\lambda_{k+1},\dots,\lambda_m]]
$$
is also dense. Therefore the functional calculus with the property stated in the theorem is unique.

\begin{proof}[Proof of Theorem~\ref{t3.3}]
By Theorem~\ref{t3.1}, the formula \eqref{eq3.2} defines an ordered $C^\infty$-functional calculus $\theta$ for~$\mathbf{b}$. Thus the proof reduces to checking that $\theta$ factors through $C^\infty(\mathbb{R}^k)\mathbin{\widehat{\otimes}} \mathbb{R}[[\lambda_{k+1},\dots,\lambda_m]]$.

Let $I$ be the closed ideal of $C^\infty(\mathbb{R}^{m-k})$ consisting of functions  all of whose partial derivatives vanishing at the origin. Identifying $C^\infty(\mathbb{R}^k)\mathbin{\widehat{\otimes}} C^\infty(\mathbb{R}^{m-k})$ with $C^\infty(\mathbb{R}^m)$, we can treat $C^\infty(\mathbb{R}^k)\mathbin{\widehat{\otimes}} I$ as a subspace of $C^\infty(\mathbb{R}^m)$. By~\eqref{eq3.4}, it suffices to show that $\theta(F)=0$ for every $F$ in $C^\infty(\mathbb{R}^k)\mathbin{\widehat{\otimes}} I$. Moreover, we can assume that $F=F_1\otimes F_2$, where $F_1\in C^\infty(\mathbb{R}^k)$ and $F_2\in I$.

Let $f_1\in \mathcal{S}(\mathbb{R}^k)$ and $f_2\in \mathcal{S}(\mathbb{R}^{m-k})$ be pre-images of $F_1$ and $F_2$. Then \eqref{eq3.2} implies that
$$
\theta(F_1\otimes F_2)=f_1(b_1,\dots,b_k)f_2(b_{k+1},\dots,b_m).
$$
We claim that the second multiple is~$0$.

Since $b_{k+1},\dots,b_m$ are nilpotent, there is $n\in\mathbb{N}$ such that
$$
b_{k+1}^{n+1}=\cdots=b_m^{n+1}=0.
$$
Write the formula in~\eqref{eq3.2} for $f_2(b_{k+1},\dots,b_m)$. Expanding each of the exponentials into Taylor series, we see that the integral on the right-hand side is equal to
\begin{equation}
\label{eq3.5}
\sum_{\alpha_{k+1},\dots, \alpha_m=0}^n \int_{\mathbb{R}^{m-k}} \widehat f_2(\mathbf{s}) \frac{(i s_{k+1} b_{k+1})^{\alpha_{k+1}}\cdots (i s_m b_m)^{\alpha_m}}{\alpha_{k+1}!\cdots \alpha_m!}\,d\mathbf{s}.
\end{equation}
Since $\widehat{f_2^{(\alpha)}}(\mathbf{s})=(i s_{k+1})^{\alpha_{k+1}}\cdots (i s_m)^{\alpha_m}\, \widehat f_2(\mathbf{s})$ for every multi-index $\alpha\in \mathbb{Z}_+^{m-k}$, the integrand in each summands has the form
$$
\widehat{f_2^{(\alpha)}}(\mathbf{s}) \frac{b_{k+1}^{\alpha_{k+1}}\cdots b_m^{\alpha_m}}{\alpha_{k+1}!\cdots \alpha_m!}.
$$
Since $F_2\in I$, all the partial derivatives of $f_2$ vanish at the origin. Taking for every $\alpha$ the inverse Fourier transform of~$f_2^{(\alpha)}$, we obtain
$$
\int_{\mathbb{R}^{m-k}} \widehat{f_2^{(\alpha)}}(\mathbf{s})\,d\mathbf{s}=0.
$$
Hence every summand in~\eqref{eq3.5} is~$0$. Thus $f_2(b_{k+1},\dots,b_m)=0$ and so $\theta(F_1\otimes F_2)=0$.
\end{proof}

\section{Globally defined non-commutative functions of class $C^\infty$ and a~multiplicative calculus for triangular finite-dimensional Lie algebras}
\label{s4}

\subsection{Triangular finite-dimensional Lie algebras}
\label{ss4.1}

Recall that a Lie algebra $\mathfrak{g}$ over a field~$\mathbb{K}$ is said to be \textit{triangular} if it is soluble and for every $x\in\mathfrak{g}$ all the eigenvalues of the linear operator~$\operatorname{ad} x$ are  in~$\mathbb{K}$\footnote{Note that  for $\mathfrak{g}$ to be triangular  it suffices to satisfy the first condition in the case when $\mathbb{K}=\mathbb{C}$ (by Lie's theorem) and the second condition in the case when  $\mathbb{K}=\mathbb{R}$ (see the discussion 'About supersolvable Lie algebras' on https:/\!/mathoverflow.net/questions/207154/about-supersolvable-lie-algebras).}.

The choice of the class of triangular Lie algebras as the main is justified by the following assertion.

\begin{proposition}
\label{p4.1}
Let $\mathfrak{g}$ be a finite-dimensional Lie $\mathbb{R}$-algebra, $B$ a Banach $\mathbb{R}$-algebra of polynomial growth, and $\gamma\colon \mathfrak{g} \to B$ a Lie $\mathbb{R}$-algebra homomorphism. Then $\gamma(\mathfrak{g})$ is a triangular Lie algebra.
\end{proposition}

\begin{proof}
Denote $\gamma(\mathfrak{g})$ by $\mathfrak{h}$. It follows form Theorem~\ref{t2.8} that the ideal $[\mathfrak{h},\mathfrak{h}]$ is contained in $\operatorname{Rad} B$. By Proposition~\ref{p2.9}, $\operatorname{Rad} B$ is nilpotent. Therefore $[\mathfrak{h},\mathfrak{h}]$ is nilpotent as a Lie algebra. Hence $\mathfrak{h}$ is soluble.

Next, we use a standard trick from the theory of generalized scalar operators. Consider an operator $\operatorname{ad} h\colon B\to B$ and note that
$$
e^{is\operatorname{ad} h}(b)=e^{ish}b e^{-ish},\qquad s\in\mathbb{R},\quad h\in \mathfrak{h},\quad b\in B;
$$
for a proof see, for example,  \cite[Chapter~II, \S\,15, p.\,83, Remark~1]{29}. Since $h\in\mathfrak{h}$ is of polynomial growth, so is the operator $\operatorname{ad} h$. Since $\mathfrak{h}$ is an invariant subspace of $\operatorname{ad} h$, the restriction to $\mathfrak{h}$ is also of polynomial growth. In particular, the spectrum of $(\operatorname{ad} h)|_\mathfrak{h}$ is in~$\mathbb{R}$, which means, by definition, that~$\mathfrak{h}$ is triangular.
\end{proof}

\begin{corollary}
\label{c4.2}
Let $\mathfrak{g}$ be a finite-dimensional Lie $\mathbb{R}$-algebra, $B$ a projective limit of Banach $\mathbb{R}$-algebras of polynomial growth, and $\gamma\colon \mathfrak{g} \to B$ a Lie $\mathbb{R}$-algebra homomorphism. If~$\gamma$ is injective, then~$\mathfrak{g}$ is triangular.
\end{corollary}

\begin{proof}
Since $\mathfrak{g}$ is finite dimensional, there is a submultiplicative seminorm $\|\,{\cdot}\,\|$ on~$B$ such that its restriction to~$\gamma(\mathfrak{g})$ is a norm. The completion of~$B$ with respect to~$\|\,{\cdot}\,\|$ is a Banach algebra, to which we can apply Proposition~\ref{p4.1}.
\end{proof}

\subsection{Statement of main results}
\label{ss4.2}
Let $\mathfrak{g}$ be a triangular finite-dimensional Lie $\mathbb{R}$-algebra. Denote by~$\mathfrak{n}$ the nilpotent radical  of~$\mathfrak{g}$ (that is, the intersection  of kernels of all finite-dimensional irreducible representations). Since $\mathfrak{g}$ is soluble,  the well-known formula $\mathfrak{n}=[\mathfrak{g},\mathfrak{g}]$ holds. Fix a linear basis $e_{k+1},\dots, e_m$ of~$\mathfrak{n}$ with some complement $e_1,\dots, e_k$ to a linear basis of~$\mathfrak{g}$ and consider the Fr\'echet space
\begin{equation}
\label{eq4.1}
C^\infty_\mathfrak{g}:=C^\infty(\mathbb{R}^k)\mathbin{\widehat{\otimes}} \mathbb{R}[[e_{k+1},\dots,e_m]].
\end{equation}
We use the corresponding PBW basis $\{e^\alpha:=e_1^{\alpha_1}\cdots e_m^{\alpha_m}\colon \alpha\in \mathbb{Z}_+^m\}$ of~$U(\mathfrak{g})$ and the linear map $\Phi$ (see \eqref{eq1.1}) to identify $U(\mathfrak{g})$ with a~subspace of~$C^\infty_\mathfrak {g}$ (which is dense by the Weierstrass approximation theorem for continuously differentiable functions). It is shown below (Theorem~\ref{t4.4}) that the space $C^\infty_\mathfrak{g}$ defined in this way is  independent of the choice of the basis.

The following two theorems are main results of the paper.

\begin{theorem}
\label{t4.3}
Let $\mathfrak{g}$ be a triangular finite-dimensional real Lie algebra. Then the multiplication in~$U(\mathfrak{g})$  extends to a continuous multiplication in~$C^\infty_\mathfrak{g}$ (for every basis in~$\mathfrak{n}$ and every complement to a basis in~$\mathfrak{g}$). Moreover,   with this multiplication, $C^\infty_\mathfrak{g}$ is a projective limit of real Banach algebras of polynomial growth and hence a Fr\'echet--Arens--Michael algebra of polynomial growth.
\end{theorem}

The restriction of the embedding $U(\mathfrak{g})\to C^\infty_\mathfrak{g}$ to $\mathfrak{g}$ is a real Lie algebra homomorphism $ \mathfrak{g}\to C^\infty_\mathfrak{g}$. Denote it by $\mu$. The following result is a generalization of Theorem~\ref{t3.2}.

\begin{theorem}
\label{t4.4}
Let $\mathfrak{g}$ be a triangular finite-dimensional real Lie algebra. If  $B$ is a projective limit of real Banach algebras of polynomial growth  and $\gamma\colon \mathfrak{g}\to B$ is a real Lie algebra homomorphism, then
there exists a multiplicative $C^\infty_\mathfrak{g}$-functional calculus, that is, a continuous homomorphism $\theta\colon C^\infty_\mathfrak{g}\to B$ such that the diagram
\begin{equation}
\label{eq4.2}
\begin{gathered}
\xymatrix@C=35pt{
\mathfrak{g}\ar[d]_\mu\ar[dr]^{\gamma} &
\\
C^\infty_\mathfrak{g}\ar[r]_{\theta} & B
}
\end{gathered}
\end{equation}
is commutative. As a corollary, the algebra $C^\infty_\mathfrak{g}$ is independent of the choice of a
basis in~$\mathfrak{n}$ and a complement to a basis in~$\mathfrak{g}$
\end{theorem}

Theorems~\ref{t4.3} and~\ref{t4.4} immediately imply the following result.

\begin{corollary}
\label{c4.5}
Every homomorphism $\mathfrak{g}\to\mathfrak{h}$ of triangular finite-dimensional real Lie algebras induces a continuous homomorphism $C^\infty_\mathfrak{g}\,{\to}\, C^\infty_\mathfrak{h}$.
The resulting correspondence is a functor from the category of triangular real Lie algebras to the category of Fr\'echet--Arens--Michael $\mathbb{R}$-algebras of polynomial growth.
\end{corollary}

Thus, all three requirements for~$C^\infty_\mathfrak{g}$ listed in the introduction are satisfied: (1)~if $\mathfrak{g}=\mathbb{R}^m$ for some $m\in\mathbb{N}$, then $C^\infty_\mathfrak{g}=C^\infty(\mathbb{R}^m)$; (2)~the algebra $C^\infty_\mathfrak{g}$ is a completion of~$U(\mathfrak{g})$ and homomorphism $U(\mathfrak{g})\to C^\infty_\mathfrak{g}$ is injective (by Corollary~\ref{c4.2}); (3)~the correspondence  $\mathfrak{g}\mapsto C^\infty_\mathfrak{g}$  extends to a functor.

The definition of the space $C^\infty(\mathbb{R}^k)\mathbin{\widehat{\otimes}} \mathbb{R}[[e_{k+1},\dots,e_m]]$ in~\eqref{eq4.1} can be applied to an arbitrary soluble algebra~$\mathfrak{g}$. However, as can be seen from the following example, in the case when $\mathfrak{g}$ is not triangular, the multiplication in~$U(\mathfrak{g})$ need not be extendable to a continuous operation on this space.

\begin{example}
\label{e4.6}
Consider the Lie algebra $\mathfrak{e}_2$ of the group of motions of the plane~$\mathbb{R}^2$. It has a linear basis $e_1,e_2,e_3$ with relations
$$
[e_1,e_2]=e_3, \quad [e_1,e_3]=-e_2, \quad\text{and}\quad [e_2,e_3]=0.
$$
This algebra is the simplest example of a soluble non-triangular Lie algebra (the eigenvalues of the operator $\operatorname{ad} e_1$ are $i$, $0$, $-i$).

By the PBW theorem, $U(\mathfrak{e}_2)$ is a subspace of $C^\infty(\mathbb{R})\mathbin{\widehat{\otimes}}\mathbb{R}[[e_2,e_3]]$. Proposition~\ref{p4.1} immediately implies that the last space  cannot be an algebra of polynomial growth with the multiplication inherited from $U(\mathfrak{e}_2)$. But to prove that it is not even a Fr\'echet algebra we need additional reasoning.
Let $f\in\mathbb{R}[\lambda]$. It follows from a well-known commutation formula  (see, for example, \cite[\S\,15, p.\,82, Corollary~1]{29}) that
$$
\operatorname{ad} f(e_1)(e_2) = \sum_{n\ge 1}\frac{(-1)^{n+1}}{n!} f^{(n)}(e_1)(\operatorname{ad} e_1)^n(e_2).
$$
Since $(\operatorname{ad} e_1)^{2k-1}(e_2)=(-1)^{k+1} e_3$ and $(\operatorname{ad} e_1)^{2k}(e_2)=(-1)^k e_2$ for every $k\in\mathbb{N}$, we have
$$
e_2f(e_1) = f(e_1)e_2+ \sum_{k\ge 1}\frac{(-1)^k}{(2k-1)!} f^{(2k-1)}(e_1)e_3+\sum_{k\ge 1}\frac{(-1)^k}{(2k)!} f^{(2k)}(e_1)e_2.
$$

Treating $f$ as a polynomial in a complex variable, we obtain from Taylor's formula that
\begin{equation}
\label{eq4.3}
e_2f(e_1) = \frac{f(e_1-i)-f(e_1+i)}{2i}\,e_3+\frac{f(e_1-i)+f(e_1+i)}{2}\,e_2.
\end{equation}
Put $f(\lambda)=(\lambda^2+1)^{-1}$ and consider the sequence of compact subsets~$K_m$ of~$\mathbb{C}$ defined by
$$
K_m:=\biggl\{\lambda\in\mathbb{C}\colon \,|\mathbb{R}e \lambda|\le m,\, |{\operatorname{Im} \lambda}|\le 1,\, |\lambda-i|\ge \frac1{m},\, |\lambda+i|\ge \frac1{m}\biggr\},\qquad m\in\mathbb{N}.
$$
Let $f_n$ be a sequence of polynomials converging to~$f$ uniformly on each of~$K_m$ (the existence follows from Mergelyan's theorem; see, for example, \cite[Theorem 20.5]{30}). Moreover, we can assume that every $f_n$ has real coefficients. Since
$$
\frac{f(\lambda-i)-f(\lambda+i)}{2i}=\frac{2}{\lambda^3+4\lambda},
$$
the sequence $f_n(\lambda-i)-f_n(\lambda+i)$ converges  uniformly in~$\lambda$ on every compact subset of~$\mathbb{R}\setminus\{0\}$ to a rational function with pole at~$0$.

Now assume the opposite, that is, that the multiplication in~$U(\mathfrak{e}_2)$  continuously extends to $C^\infty(\mathbb{R})\mathbin{\widehat{\otimes}}\mathbb{R}[[e_2,e_3]]$.
Then $e_2f_n(e_1)$ converges to~$e_2f(e_1)$ in the topology inherited from $C^\infty(\mathbb{R})\mathbin{\widehat{\otimes}}\mathbb{R}[[e_2,e_3]]$.  We have from~\eqref{eq4.3} that $f_n(\lambda-i)-f_n(\lambda+i)$ converges  uniformly in~$\lambda$ on  compact subsets to a~continuous function on~$\mathbb{R}$. We get a contradiction, which means that there is no continuous extension.
\end{example}

Nevertheless, it is possible to give a more general definition of $C^\infty_\mathfrak{g}$, which works for an arbitrary Lie $\mathbb{R}$-algebra $\mathfrak{g}$ (not only soluble and not only finite-dimensional); see~\cite{5}. Of course, when non-triangular $\mathfrak{g}$, it does not coincide with~\eqref{eq4.1} and $U(\mathfrak{g})\to C^\infty_\mathfrak{g}$ is not injective.

In the rest of this section we are mainly concerned with proving Theorems~\ref{t4.3} and~\ref{t4.4}. The proof of the first of them is contained in \S\S\,\ref{ss4.3}--\ref{ss4.6} and the proof of the second in \S\,\ref{ss4.7}.

\subsection{Beginning of the proof of Theorem~\ref{t4.3}}
\label{ss4.3}
We start to prove Theorem~\ref{t4.3}. In this section, we describe preliminary constructions and reduce the assertion of the theorem to verifying the topological injectivity of a certain homomorphism~$\rho$. Further, in \S\,\ref{ss4.4} we consider two examples and give proofs in those cases as an illustration of essential features of our argument. Auxiliary assertions needed in the general case are contained in~\S\,\ref{ss4.5}. The proof of the topological injectivity of~$\rho$ is completed in~\S\,\ref{ss4.6}.

Let $e_1,\dots, e_m$ be a linear basis of~$\mathfrak{g}$ such that $e_{k+1},\dots, e_m$ is a linear basis of~$\mathfrak{n}$. We endow $U(\mathfrak{g})$ with the topology inherited from $C^\infty_\mathfrak{g}$; see~\eqref{eq4.1}. The idea of the proof is to construct a \text{topologically injective}
homomorphism from $U(\mathfrak{g})$ to a product of Banach $\mathbb{R}$-algebras of polynomial growth. This is sufficient since the property of having polynomial growth is preserved  under passage to a closed subalgebra (see Proposition~\ref{p2.11}).

We first check that the topology on $C^\infty_\mathfrak{g}$ given above is independent of the choice\footnote{It also is independent of the choice of a complement to a basis of~$\mathfrak{g}$ but this is shown below in the proof of the Theorem~\ref{t4.4}. The current argument can be applied to any choice of a complement and the independence is not used here.}
of a basis of~$\mathfrak{n}$. Let $e'_{k+1},\dots, e'_m$ be another basis of~$\mathfrak{n}$. Consider two topologies on  $U(\mathfrak{g})$ given by the embeddings
$$
C^\infty(\mathbb{R}^k)\mathbin{\widehat{\otimes}} \mathbb{R}[[e_{k+1},\dots,e_m]]\quad\text{and}\quad C^\infty(\mathbb{R}^k)\mathbin{\widehat{\otimes}} \mathbb{R}[[e'_{k+1},\dots,e'_m]],
$$
respectively. We need to show that these topologies coincide.

Every $a\in U(\mathfrak{g})$ can be expanded into a sum in two ways  by those elements of the PBW basis that do not contain $e_1,\dots, e_k$  with  non-commutative polynomials as coefficients, namely:
\begin{equation}
\label{eq4.4}
a=\sum_{\alpha\in\mathbb{Z}_+^{m-k}} \Phi(f_\alpha)e^\alpha =\sum_{\beta\in\mathbb{Z}_+^{m-k}} \Phi(h_\beta)(e')^\beta
\end{equation}
(here we identify  $\mathbb{Z}_+^{m-k}$ with the set of multi-indices $\alpha\in\mathbb{Z}_+^{m}$ such that $\alpha_1=\dots=\alpha_k=0$), where $f_\alpha,h_\beta\in \mathbb{R}[\lambda_1,\dots,\lambda_k]$ and $\Phi\colon \mathbb{R}[\lambda_1,\dots,\lambda_k]\to U(\mathfrak{g})$ is an ordered calculus defined as in~\eqref{eq1.1}.

As above, in the commutative case (that is, when $p=1$) we denote by $|\,{\cdot}\,|_{K,n}$ the seminorms $\|\,{\cdot}\,\|_{1,K,n}$ defined in~\eqref{eq2.2}. It is easy to see that the topology on $U(\mathfrak{g})$ corresponding to $e_1,\dots, e_m$ is determined by the family
$$
|a|_{\alpha,K,n}:= |f_\alpha|_{K,n} \qquad (n\in\mathbb{Z}_+,\quad \alpha \in \mathbb{Z}_+^{m-k},\quad K\subset \mathbb{R}^m\text{ and is compact})
$$
and the topology corresponding to  $e_1,\dots, e_k,e'_{k+1},\dots, e'_m$ is determined  by the family
$$
|a|'_{\beta,K,n}:= |h_\beta|_{K,n},\qquad n\in\mathbb{Z}_+,\quad \beta \in \mathbb{Z}_+^{m-k},\quad K\subset \mathbb{R}^m\text{ and is compact}.
$$
Since every monomial $(e')^\beta$ is a linear combination of $e^\alpha$, it follows from~\eqref{eq4.4} that every $\Phi(f_\alpha)$ is a linear combination of $\Phi(h_\beta)$ with coefficients independent of~$f_\alpha$. This implies that $(|\,{\cdot}\,|_{\alpha,K,n})$ is majorized by  $(|\,{\cdot}\,|'_{\beta,K, n})$. A similar argument shows that the second family is majorized by the first. Thus the topologies determined by them coincide.

Later on $e_{k+1},\dots, e_m$ will be chosen in a special way but right here we consider the following construction applicable to any such basis (cf. special cases below in~\eqref{eq4.8} and~\eqref{eq4.16}).

Given a finite-dimensional representation $\pi\colon U(\mathfrak{g})\to \operatorname{End} L$, where $L$ is a vector $\mathbb{R}$-space, we define a homomorphism
$\widetilde\pi\colon U(\mathfrak{g})\to\mathbb{R}[\lambda_1,\dots,\lambda_k]\otimes \operatorname{End} L$ by putting
\begin{equation}
\label{eq4.5}
\widetilde\pi(e_j):=\lambda_j\otimes 1+1\otimes \pi(e_j),\quad j\le k,\qquad \widetilde\pi(e_j):=1\otimes\pi(e_j), \quad j> k.
\end{equation}
(The homomorphism is well defined because $\mathfrak{n}=[\mathfrak{g},\mathfrak{g}]$.) Extending the codomain of~$\widetilde\pi$ we can assume that it coincides with $C^\infty(\mathbb{R}^k)\otimes \operatorname{End} L$. On the other hand, since every finite-dimensional representation of a triangular Lie algebra is equivalent to an upper triangular representation (by a generalization of Lie's theorem; see   \cite[\S\,1.2, p.\,11, Theorem~1.2]{31}), the range of $\widetilde\pi$ is in fact in $C^{\infty}(\mathbb{R}^k)\otimes \mathrm{T}_d$, where $d$ is the dimension of~$\pi$. We endow the last algebra with the topology of projective tensor product and so identify it with $C^{\infty}(\mathbb{R}^k)\mathbin{\widehat{\otimes}} \mathrm{T}_d$ (or $C^{\infty}(\mathbb{R}^k,\mathrm{T}_d)$, this follows from the facts that it is isomorphic to the injective tensor product \cite[Theorem 44.1]{32} and that the factors are nuclear).

It follows from Theorem~\ref{t3.3} that $\widetilde\pi$ extends to a continuous linear map  $C^\infty_\mathfrak{g}\to C^{\infty}(\mathbb{R}^k,\mathrm{T}_d)$.

Assume now that we have constructed a countable family $(\pi_\beta\colon U(\mathfrak{g})\to \operatorname{End} L_\beta)$ of finite-dimensional representations and then consider a homomorphism
\begin{equation}
\label{eq4.6}
\rho\colon U(\mathfrak{g})\to \prod_{\beta} C^{\infty}(\mathbb{R}^k,\mathrm{T}_{d_\beta}) \colon a\mapsto (\widetilde\pi_\beta(a)),
\end{equation}
where $d_\beta$ is the dimension of~$\pi_\beta$. Since every $\widetilde\pi_\beta$ is continuous,  so is $\rho$.

To complete the argument it suffices to take  a basis of~$\mathfrak{g}$ and a family $(\pi_\beta)$ of representations  such that the homomorphism $\rho$ defined in~\eqref{eq4.6} is topologically injective. Indeed, then we can identify $C^\infty_\mathfrak{g}$ with a closed subalgebra of $D:=\prod C^{\infty}(\mathbb{R}^k,\mathrm{T}_{d_\beta})$, which is an Arens-Michael-Fr\'echet algebra, and so $C^\infty_\mathfrak{g}$ is of this type. Moreover, Theorem~\ref{t2.12} implies that each factor in~$D$ is a projective limit of Banach algebras of polynomial growth and hence, by Proposition~\ref{p2.11}, it is isomorphic to a closed subalgebra of a product of such algebras. Therefore $D$ and its subalgebra $C^\infty_\mathfrak{g}$ are of this form. Applying again Proposition~\ref{p2.11}, we conclude that $C^\infty_\mathfrak{g}$ is a projective limit of Banach algebras of polynomial growth and hence, by Proposition~\ref{p2.11}, it is isomorphic to a closed subalgebra of a product of algebras of polynomial growth and hence is itself of polynomial growth.

The construction used in the proof of the general case is rather cumbersome. So we first consider two simplest examples illustrating the main ideas of the argument.

\subsection{Examples for the proof of Theorem~\ref{t4.3}}
\label{ss4.4}

\begin{example}
\label{e4.7}
Denote by~$\mathfrak{a}\mathfrak{f}_1$ the two-dimensional real Lie algebra with linear basis $e_1,e_2$ and multiplication determined by the relation $[e_1,e_2]= e_2$.
Consider the following sequence $(\pi_q,\,q\in\mathbb{Z}_+)$ of representations of $\mathfrak{a}\mathfrak{f}_1$. Take  the trivial representation as $\pi_0$. For $q\in\mathbb{N}$ put $\pi_q(e_1)=X_q$ and $\pi_q(e_2)=Y_q$, where
\begin{equation}
\label{eq4.7}
X_q:= \begin{pmatrix}
q& &&& \\
 &q-1 && &\\
 & & \ddots&&\\
 & & &1&\\
 & & &&0
\end{pmatrix},\qquad
Y_q:= \begin{pmatrix}
0& 1&&& \\
&0& 1&& \\
 && \ddots& \ddots&\\
 & &&0 &1\\
 & &&&0
\end{pmatrix}.
\end{equation}
(We identify an operator with its matrix.) Here and below, we mainly omit diagonals consisting of zeros.

Let $\Phi\colon \mathbb{R}[\lambda]\to U(\mathfrak{a}\mathfrak{f}_1)$ denote the functional calculus for $e_1$. Write an element of  $U(\mathfrak{a}\mathfrak{f}_1)$ in the form $a=\sum_j \Phi(f_j)e_2^j$, where $f_j\in \mathbb{R}[\lambda]$. Then, as is easy to check,
$$
\pi_q(a)=\begin{pmatrix}
f_0(q)&f_1(q) &\cdots&&f_q(q) \\
 &f_0(q-1) && &\\
 & & \ddots&&\vdots\\
 & & &f_0(1)&f_1(1)\\
 & & &&f_0(0)
\end{pmatrix}.
$$
The homomorphism $\widetilde\pi_q\colon U(\mathfrak{a}\mathfrak{f}_1)\to\mathbb{R}[\lambda]\otimes \mathrm{T}_{q+1}$ defined in~\eqref{eq4.5} take the form
\begin{equation}
\label{eq4.8}
\widetilde\pi_q(e_1):=\lambda\otimes 1+1\otimes \pi_q(e_1),\qquad \widetilde\pi_q(e_2):=1\otimes\pi_q(e_2) .
\end{equation}
Identifying $\mathbb{R}[\lambda]\otimes \mathrm{T}_{q+1}$ with the algebra $\mathbb{R}[\lambda;\mathrm{T}_{q+1}]$  of matrix-valued polynomial functions, we can write
\begin{equation}
\label{eq4.9}
[\widetilde\pi_q(a)](\lambda)=\begin{pmatrix}
f_0(\lambda+q)&f_1(\lambda+q) &\cdots&&f_q(\lambda+q) \\
 &f_0(\lambda+q-1) && &\\
 & & \ddots&&\vdots\\
 & & &f_0(\lambda+1)&f_1(\lambda+1)\\
 & & &&f_0(\lambda)
\end{pmatrix}.
\end{equation}

It follows from Theorems~\ref{t2.12} and~\ref{t3.3} that the functional calculus $C^\infty(\mathbb{R})\mathbin{\widehat{\otimes}} \mathbb{R}[[e_2]]\to C^\infty(\mathbb{R},\mathrm{T}_{q+1})$ is continuous for every~$q$. Consider the continuous homomorphism defined by
\begin{equation}
\label{eq4.10}
\rho\colon U(\mathfrak{a}\mathfrak{f}_1)\to \prod_{q=0}^\infty C^\infty(\mathbb{R},\mathrm{T}_{q+1})
\end{equation}
(cf.~\eqref{eq4.6}). We claim that it is topologically injective.

The topology on $U(\mathfrak{a}\mathfrak{f}_1)$ inherited from $C^\infty(\mathbb{R})\mathbin{\widehat{\otimes}} \mathbb{R}[[e_2]]$ can be determined by the following family of seminorms:
\begin{equation}
\label{eq4.11}
|a|_{q,M,l}:= |f_q|_{M,l},\qquad q,l \in \mathbb{Z}_+,\quad M\subset \mathbb{R}\text{ and is compact}.
\end{equation}
Moreover, it is sufficient to take as $M$ only intervals, $[c,d]\subset \mathbb{R}$. On the other hand, for $n,p \in \mathbb{Z}_+$ and $[\gamma,\delta]\subset \mathbb{R}$  the seminorm $\|\widetilde\pi_p(\,{\cdot}\,)\|_{p+1,[\gamma,\delta],n}$ on $U(\mathfrak{a}\mathfrak{f}_1)$, where $\|\,{\cdot}\,\|_{p+1,[\gamma,\delta],n}$ is defined in~\eqref{eq2.2},  is clearly continuous with respect to the topology  inherited from the codomain of~$\rho$ (here we take the operator norm as $\|\,{\cdot}\,\|_{p+1}$).

Since $f_q$ corresponds to the upper-right entry  in~\eqref{eq4.9}, it follows that
\begin{equation}
\label{eq4.12}
|f_q|_{[c,d],l}\le\|\widetilde\pi_q(a) \|_{q+1,[c-q,d-q],l}.
\end{equation}
Thus, given $q$, $[c,d]$ and $l$, the seminorm $|\,{\cdot}\,|_{q,[c,d],l}$ is majorized  (in the sense of \cite[Definition~4.1.2]{33}) by the family $\{\|\widetilde\pi_p(\,{\cdot}\,)\|_{p+1,[\gamma,\delta],n}\}$. This means that $\rho$ is topologically injective and the assertion of Theorem~\ref{t4.3} is proved for~$\mathfrak{a}\mathfrak{f}_1$.
\end{example}

\begin{example}
\label{e4.8}
Consider the three-dimensional real Heisenberg algebra, that is, the Lie algebra $\mathfrak{h}$ with linear basis  $e_1,e_2,e_3$ and multiplication determined by the relations
$[e_1,e_2]=e_3$, $[e_1,e_3]=[e_2,e_3]=0$.

Denote by $\pi_0$ the trivial representation of $\mathfrak{h}$, by  $\pi_1$ the
`standard representation',
\begin{equation}
\label{eq4.13}
\pi_1(e_1):=
\begin{pmatrix} 0 & 1 & 0
\\ 0 & 0 & 0 \\
0&0& 0\end{pmatrix}, \qquad \pi_1(e_2):=
\begin{pmatrix} 0 & 0 & 0
\\ 0 & 0 & 1 \\
0&0&0 \end{pmatrix}, \qquad
 \pi_1(e_3):=
\begin{pmatrix} 0 & 0 &1
\\ 0 & 0 & 0\\
0&0&0 \end{pmatrix},
\end{equation}
and by $\pi_q$ the $q$-th tensor power of~$\pi_1$.

In contrast to Example~\ref{e4.7},  to prove the topological injectivity of the homomorphism
$$
\rho\colon U(\mathfrak{h})\to \prod_{q=0}^\infty C^\infty(\mathbb{R}^2,\mathrm{T}_{q})
$$
we need an argument by induction.

The topology on $U(\mathfrak{h})$ inherited from $C^\infty_\mathfrak{h}$ can be determined by the following family of seminorms:
\begin{equation}
\label{eq4.14}
|a|_{q,M,l}:= |f_q|_{M,l},\qquad q,l\in\mathbb{Z}_+,\quad M\subset \mathbb{R}^2\text{ and is compact}.
\end{equation}
It suffices to show that, given $q$, $M$ and $l$, the seminorm $|\,{\cdot}\,|_{q,M,l}$ can be majorized by the family
$(\|\widetilde\pi_p(\,{\cdot}\,)\|_{p+1,K,n})$, where $\|\,{\cdot}\,\|_{p+1,K,n}$ is defined in~\eqref{eq2.2} (here we again take the operator norm as  $\|\,{\cdot}\,\|_{p+1}$).

Let $\Phi\colon \mathbb{R}[\lambda,\mu]\to U(\mathfrak{h})$ denote the ordered functional calculus for  $e_1$ and $e_2$. Write an element $U(\mathfrak{h})$ in the form $a=\sum_j \Phi(f_j)e_3^j$, where $f_j\in \mathbb{R}[\lambda,\mu]$.

We proceed by induction on $q$. When $q=0$, note that
$$
\widetilde\pi_0\colon U(\mathfrak{h})\to\mathbb{R}[\lambda,\mu]\colon e_1\mapsto\lambda,\,e_2\mapsto\mu,\, e_3\mapsto 0.
$$
Then $\widetilde\pi_0(a)=f_0$ and it is clear that
\begin{equation}
\label{eq4.15}
|a|_{0,M,l}=\|\widetilde\pi_0(a)\|_{1,M,l}.
\end{equation}
Thus the assertion holds for $q=0$.

When $q=1$, we have that $\widetilde\pi_1\colon U(\mathfrak{h})\to\mathbb{R}[\lambda,\mu]\otimes \mathrm{T}_3$ in~\eqref{eq4.5} has the form
\begin{equation}
\label{eq4.16}
\begin{aligned}
\widetilde\pi_1(e_1) &=\lambda\otimes 1+1\otimes \pi_1(e_1),
\\
\widetilde\pi_1(e_2) &=\mu\otimes 1+1\otimes \pi_1(e_2),
\\
\widetilde\pi_1(e_3) &=1\otimes\pi_1(e_3).
\end{aligned}
\end{equation}
It is easy to see that
$$
\widetilde\pi_1(a)=\begin{pmatrix}
f_0 & \dfrac{\partial f_0}{\partial \lambda } & \dfrac{\partial^2 f_0}{\partial\lambda \,\partial\mu}+f_1
\\[3mm]
0 & f_0 & \dfrac{\partial f_0}{\partial \mu}
\\[2mm]
0&0&f_0 \end{pmatrix}.
$$

It follows from~\eqref{eq4.14} that $|a|_{1,M,l}=|f_1|_{M,l}$. To estimate $|f_1|_{M,l}$ consider the linear functional $\eta$  on $\mathrm{T}_3$ that maps a matrix to its upper-right entry. Then
$$
f_1=(1\otimes \eta)(\widetilde\pi_1(a)) -\frac{\partial^2 f_0}{\partial\lambda \,\partial\mu}.
$$
It is clear that
$|(1\otimes \eta)(\widetilde\pi_1(a))|_{M,l}\le \|\widetilde\pi_1(a)\|_{3,M,l}$ and
$$
\biggl|\frac{\partial^2 f_0}{\partial\lambda \,\partial\mu}\biggr|_{M,l} =|f_0|_{M,l+2}=\|\widetilde\pi_0(a)\|_{1,M,l+2}.
$$
Therefore,
\begin{equation}
\label{eq4.17}
|a|_{1,M,l}=|f_1|_{M,l}\le \|\widetilde\pi_1(a)\|_{3,M,l}+\|\widetilde\pi_0(a)\|_{1,M,l+2}.
\end{equation}
Thus $|\,{\cdot}\,|_{1,M,l}$ is majorized by the family of seminorms defined above. This completes the inductive step  from $q=0$ to~$q=1$. (The bound for the second summand in~\eqref{eq4.17} can also be obtained from general considerations involving the continuity of $\widetilde\pi_1$  with respect to the topology on $C^\infty_\mathfrak{h}$, which implied by Theorem~\ref{t3.3}. This is how we reason in the general case.)

Now assume that $q>1$ and for every $M$, $l$  the seminorm $|\,{\cdot}\,|_{q',M,l}$ is majorized  by the family $(\|\widetilde\pi_p(\,{\cdot}\,)\|_{p+1,K,n})$ when $q'<q$. We need to show that
for every  $M$ and $l$ the seminorm $|\,{\cdot}\,|_{q,M,l}$ is also majorized  by this family.

Since $\pi_q$ is the $q$-th tensor power of~$\pi_1$, we have
$$
\pi_q(e_j)=\pi_1(e_j)\otimes 1\otimes\dots\otimes 1+\cdots+ 1\otimes\dots\otimes 1 \otimes \pi_1(e_j),
$$
where the sum contains $q$ terms. It follows that
\begin{equation}
\label{eq4.18}
\pi_q(e_3^q)=q! \,\pi_1(e_3)\otimes\dots\otimes\pi_1(e_3).
\end{equation}
Besides, since
\begin{equation}
\label{eq4.19}
\pi_1(e_1)\pi_1(e_3)=\pi_1(e_2)\pi_1(e_3)=0,
\end{equation}
we conclude that $\pi_q(e_1)\pi_q(e_3^q)=\pi_q(e_2)\pi_q(e_3^q)=0$. It easily follows that
$$
\widetilde\pi_q(\Phi(f)e_3^q)= f(\lambda,\mu) \otimes \pi_q(e_3^q)
$$
for every $f\in \mathbb{R}[\lambda,\mu]$ and then, noting that  $\pi_q$ is of  dimension $3^q$, we obtain
\begin{equation}
\label{eq4.20}
\|\widetilde\pi_q(\Phi(f)e_3^q)\|_{3^q,M,l}=q!\,|f|_{M,l}
\end{equation}
because $\|\pi_q(e_3^q)\|_{3^q}=q!$.

Consider the  projection
$$
P_q\biggl(\sum_n \Phi(f_n)e_3^n\biggr):=\sum_{n=0}^{q-1} \Phi(f_n)e_3^n
$$
on $C^\infty_\mathfrak{h}$.
Since $\pi_1(e_3)^2=0$, we have $\pi_q(e_3^{q'})=0$ when $q'>q$ and so
$$
\widetilde\pi_q(\Phi(f_q)e_3^q)=\widetilde\pi_q(a)-\widetilde\pi_q P_q(a)
$$
for every $a\in C^\infty_\mathfrak{h}$. Thus, the above equation and~\eqref{eq4.20} yield that
\begin{align*}
|a|_{q,M,l} &=|f_q|_{M,l}=(q!)^{-1}\|\widetilde\pi_q(\Phi(f_q)e_3^q)\|_{3^q,M,l}
\\
&\le  (q!)^{-1}(\|\widetilde\pi_q(a)\|_{3^q,M,l}+\|\widetilde\pi_q P_q(a)\|_{3^q,M,l}).
\end{align*}

Since $\widetilde\pi_q$ is continuous, $\|\widetilde\pi_q P_q(\,{\cdot}\,)\|_{3^q,M,l}$ is majorized by the family $(|P_q(\,{\cdot}\,)|_{q',M',l'})$. Since, by the induction hypothesis,  $|P_q(a)|_{q',M',l'}=|a|_{q',M',l'}$ when $q'<q$ and $|P_q(a)|_{q',M',l'}=0$ when $q'\ge q$ for every $a\in C^\infty_\mathfrak{h}$, the family  $\|\widetilde\pi_q P_q(\,{\cdot}\,)\|_{3^q,M,l}$ is also majorized by $(\|\widetilde\pi_p(\,{\cdot}\,)\|_{p+1,K,n})$. It follows that $|\,{\cdot}\,|_{q,M,l}$ is majorized by the last family and this completes the inductive step  from  $q-1$ to~$q$. (Note that in the above proof, in the case $q=1$, we can also use the projection $P_1$ instead of the functional $\eta$ and then the argument becomes a special case of the general one.)

Thus $\rho$  topologically injective and the assertion of Theorem~\ref{t4.3} holds  for~$\mathfrak{h}$.
\end{example}


\subsection{Results auxiliary to the proof of Theorem~\ref{t4.3}}
\label{ss4.5}

To implement in general the idea of the proof outlined in Examples~\ref{e4.7} and~\ref{e4.8} we need a number of auxiliary assertions. Since the complex case is also of interest, in what follows we consider a Lie algebra over a field $\mathbb{K}$ that can be either~$\mathbb{R}$ or~$\mathbb{C}$.

Let $\mathfrak{g}$ be a triangular Lie algebra and $\mathfrak{n}$ denote its nilpotent radical. We need three conditions on $x\in \mathfrak{g}$ and a finite-dimensional representation~$\pi$ of~$\mathfrak{g}$:

(A1)~$\pi(x)\ne 0$;

(A2)~for every $y\in\mathfrak{g}$ there is $\mu_y\in \mathbb{K}$ such that $\pi(y)\pi(x)=\mu_y \pi(x)$;

(A3)~for every $y\in\mathfrak{n}$ the equality $\pi(y)\pi(x)= 0$ holds, that is, $\mu_y=0$.

If the representation $\pi$ is equivalent to an upper-triangular one and is indecomposable, then it is easy to find $x$ with such properties. However, below we need  representations and vectors not only satisfying (A1)--(A3) but also additional relations, in particular, the property that the vectors form a basis of the nilpotent radical. Namely, the following generalization of~\eqref{eq4.19}  holds. The proof requires additional efforts and is connected with a modification of Ado's theorem.

\begin{proposition}
\label{p4.9}
Let $\mathfrak{g}$ be a triangular Lie algebra over~$\mathbb{K}$. Then there are a linear basis $e_{k+1},\dots, e_m$ in~$\mathfrak{n}$ and a tuple  $\pi_{k+1},\dots,\pi_m$ of finite-dimensional representations of~$\mathfrak{g}$ such that $e_r$ and~$\pi_r$ satisfies Conditions {\rm(A1)}--{\rm(A3)} for every $r\in \{k+1,\dots, m\}$ and, moreover, $\pi_r(e_j)=0$ when $r<j$.
\end{proposition}

For the proof, we need several lemmas.

\begin{lemma}
\label{l4.10}
Let $\mathfrak{g}$ be a finite-dimensional Lie algebra over~$\mathbb{K}$ and $\mathfrak{h}$ its soluble ideal.

{\rm(A)}~Then $\operatorname{rad} (\mathfrak{g}/\mathfrak{h})$ \emph{(}the soluble radical of $\mathfrak{g}/\mathfrak{h}$\emph{)} equals~$\mathfrak{r}/\mathfrak{h}$.

{\rm(B)}~If the soluble radical of~$\mathfrak{g}$ is triangular, then so is $\operatorname{rad} (\mathfrak{g}/\mathfrak{h})$.
\end{lemma}

\begin{proof}
(A)~Since $\mathfrak{h}$ is soluble, $\mathfrak{h}\subset \mathfrak{r}$. Then $\mathfrak{r}/\mathfrak{h}$  is soluble and so is contained in $\operatorname{rad} (\mathfrak{g}/\mathfrak{h})$. To check the reverse inclusion put $\mathfrak{m}:=\sigma^{-1}(\operatorname{rad} (\mathfrak{g}/\mathfrak{h}))$, where $\sigma$ denote the projection $\mathfrak{g}\to\mathfrak{g}/\mathfrak{h}$. Let $\mathfrak{s}$ be a Levi complement. It clear that $\mathfrak{r}\subset \mathfrak{m}$ and so $\mathfrak{g}=\mathfrak{r}+\mathfrak{s}$ implies that $\mathfrak{m}=\mathfrak{r}+\mathfrak{m}\cap \mathfrak{s}$. On the other hand, since $\mathfrak{s}$ is semisimple,  $\mathfrak{m}\cap \mathfrak{s}$ is a semisimple ideal of~$\mathfrak{s}$ and therefore $\sigma(\mathfrak{m}\cap \mathfrak{s})$ a semisimple subalgebra of~$\mathfrak{g}$ contained in the radical. Hence $\sigma(\mathfrak{m}\cap \mathfrak{s})=0$. Thus $\operatorname{rad} (\mathfrak{g}/\mathfrak{h})=\sigma(\mathfrak{m})=\sigma(\mathfrak{r})=\mathfrak{r}/\mathfrak{h}$.

(B)~Since $\mathbb{K}$ is of characteristic~$0$, a finite-dimensional Lie  $\mathbb{K}$-algebra is triangular if and only if when it is supersoluble, that is, contains a maximal flag of ideals (see, for example, \cite[\S\,1.2, p.\,11, Theorem~1.2]{31}). It is not hard to see that a quotient of a supersoluble Lie algebra is supersoluble. By Part~(A), this implies that  $\operatorname{rad} (\mathfrak{g}/\mathfrak{h})$ is triangular.
\end{proof}

Let $\mathfrak{z}(\mathfrak{g})$  denote the centre of $\mathfrak{g}$ and $(\mathfrak{n}_j)$ the lower central series of~$\mathfrak{n}$.

\begin{lemma}
\label{l4.11}
Let $\mathfrak{g}$ be a soluble finite-dimensional Lie algebra over~$\mathbb{K}$ and $p$ denote the nilpotency degree of $\mathfrak{n}$. Then for every non-zero $x\in \mathfrak{z}(\mathfrak{g})\cap \mathfrak{n}_{p-1}$ there is a finite-dimensional representation~$\pi$ of~$\mathfrak{g}$ such that {\rm(A1)}--{\rm(A3)} hold. If, in addition, $\mathfrak{g}$ is nilpotent, then $\pi(\mathfrak{g})$ consists of nilpotent operators.
\end{lemma}

The assertion in the case when $\mathfrak{g}$ is nilpotent is used in the proof of Proposition~\ref{p4.14}.

\begin{proof}[Proof of Lemma~\ref{l4.11}]
We construct $\pi$ for which (A1) is satisfied, as well as a condition stronger than~(A2) and~(A3), namely, $\pi(y)\pi(x)=0$ for every $y\in\mathfrak{g}$.

Our construction of the representation is the same as in the proof of Ado's theorem given in~\cite[\S\,7.4, pp.\,192--193]{34} but has a simplified form suitable for soluble Lie algebras. Consider the following abelian Lie subalgebra $\mathfrak{d}$ of the Lie algebra of derivations  of~$\mathfrak{g}$. Take the restriction of the adjoint representation of~$\mathfrak{g}$ to some Cartan subalgebra of~$\mathfrak{r}$. Then, by definition, $\mathfrak{d}$ consists of the semisimple summands of the Jordan decompositions of inner derivations associated with elements of the  Cartan subalgebra. It is clear that the restriction of an inner derivation to~$\mathfrak{z}(\mathfrak{g})$ is trivial. Since the Jordan decomposition of zero is the sum of two zeros,  the action of~$\mathfrak{d}$ on~$\mathfrak{z}(\mathfrak{g})$ is also trivial.

Put $\widehat{\mathfrak{g}}:=\mathfrak{g}\rtimes \mathfrak{d}$. Then we have a decomposition $\widehat{\mathfrak{g}}=\widehat{\mathfrak{n}}\rtimes \mathfrak{l}$, where $\widehat{\mathfrak{n}}$ is the maximal nilpotent ideal of~$\widehat{\mathfrak{g}}$ and $\mathfrak{l}$ is a reductive subalgebra \cite[\S\,7.4, p.\,193]{34}. (In the case under discussion, since $\mathfrak{g}$ is soluble, $\mathfrak{l}=\mathfrak{d}$, but we keep the notation in~\cite{34}.)

The construction of the representation $\pi$ of $\widehat{\mathfrak{g}}$ is as follows (see \cite[Proposition 7.4.4]{34}). Let $U_j$ be the two-sided ideal of~$U(\widehat{\mathfrak{n}})$ generated by the elements of the form  $n_1\cdots n_j$, where $n_1,\dots, n_j\in \widehat{\mathfrak{n}}$. Let $p$ denote the degree of nilpotency of~$\widehat{\mathfrak{n}}$. Take $U(\widehat{\mathfrak{n}})/U_p$ as the representation space and define the representation  by the formula
$$
\pi(n,d)(m+U_p):=nm+\gamma(d)(m)+U_p,\qquad n\in \widehat{\mathfrak{n}},\quad d\in \mathfrak{l},\quad m\in U(\widehat{\mathfrak{n}}),
$$
where $\gamma(d)$ is a derivation of $U(\widehat{\mathfrak{n}})$ induced by the action of~$d\in\mathfrak{l}$ on~$\widehat{\mathfrak{n}}$.

Since $\mathfrak{n}$ coincides with $[\mathfrak{g},\mathfrak{g}]$, it is invariant under the action of every derivation. Hence~$\mathfrak{n}$ is  ideal (evidently, nilpotent) of~$\widehat{\mathfrak{g}}$. Since  $\widehat{\mathfrak{n}}$ is the maximal nilpotent ideal, we have that $\mathfrak{n}\subset\widehat{\mathfrak{n}}$. In particular, $x\in \widehat{\mathfrak{n}}$. Now we show that~$x$ and the restriction of~$\pi$ to~$\mathfrak{g}$ satisfy the desired conditions. It is clear that $\pi(x,0)(m+U_p)=xm+U_p$. In particular,  $x\ne 0$ implies that $\pi(x,0)(1+U_p)\ne 0$ and so $\pi(x,0)\ne0$, that is, (A1) holds.

Now take $y=(n,d)$ in $\widehat{\mathfrak{g}}$. Then
\begin{equation}
\label{eq4.21}
\pi(n,d)\pi(x,0)(m+U_p)=nxm+\gamma(d)(xm)+U_p,\qquad m\in U(\widehat{\mathfrak{n}}).
\end{equation}

By the hypotheses of the lemma, $x\,{\in}\,\mathfrak{n}_{p-1}$ and therefore $x\,{\in}\, U_{p-1}$.
Since $U_1U_{p-1}\,{\subset}\,U_p$, we conclude
that $mnx\in U_p$ for every $m$.

Further, the action of~$\mathfrak{d}$ on~$\mathfrak{z}(\mathfrak{g})$ is trivial and so $\mathfrak{z}(\mathfrak{g})$ is in the centre of~$\widehat{\mathfrak{g}}$. In particular, the action of $\mathfrak{l}$ on~$\mathfrak{z}(\mathfrak{g})$ is also trivial. Since $x\in \mathfrak{z}(\mathfrak{g})$, we obtain $\gamma(d)(x)=0$.  Note that $\gamma(d)(m)\in U_1$ and $\gamma(d)$ is a derivation. Then it follows from  $U_{p-1}U_1\,{\subset}\, U_p$  that
$$
\gamma(d)(xm)=x\gamma(d)(m)+\gamma(d)(x)m \in U_p
$$
for every $m$. Thus~\eqref{eq4.21}  implies that $\pi(n,d)\pi(x,0)=0$. Since $\mathfrak{g}\subset \widehat{\mathfrak{g}}$, we have $\pi(y)\pi(x)=0$ for every $y\in\mathfrak{g}$.

Suppose now that  $\mathfrak{g}$ is nilpotent. Hence $\mathfrak{g}\subset \widehat{\mathfrak{n}}$ and, by the construction of~$\pi$, the set $\pi(\mathfrak{g})$ consists of nilpotent operators.
\end{proof}

\begin{lemma}
\label{l4.12}
Let $\mathfrak{g}$ be a soluble finite-dimensional Lie algebra over~$\mathbb{K}$ and $\mathfrak{h}$ an ideal of~$\mathfrak{g}$. If $x\in\mathfrak{g}/\mathfrak{h}$ and a finite-dimensional representation~$\pi$ of~$\mathfrak{g}/\mathfrak{h}$ satisfy  \emph{(A1)}--\emph{(A3)}, then every pre-image of $x$ under the projection $\sigma\colon\mathfrak{g}\to\mathfrak{g}/\mathfrak{h}$ and a representation $\pi\sigma$ also satisfy {\rm(A1)}--{\rm(A3)}.
\end{lemma}

\begin{proof}
It suffice to note that~$\mathfrak{g}/\mathfrak{h}$ is soluble and $[\mathfrak{g},\mathfrak{g}]$ maps to~$[\mathfrak{g}/\mathfrak{h},\mathfrak{g}/\mathfrak{h}]$.
\end{proof}

\begin{lemma}
\label{l4.13}
Let~$\mathfrak{g}$ be an arbitrary finite-dimensional Lie $\mathbb{K}$-algebra, and $\mathfrak{h}$ a soluble ideal of~$\mathfrak{g}$. Let $\mathfrak{k}$ denote the nilpotent radical of~$\mathfrak{g}/\mathfrak{h}$ and  $(\mathfrak{k}_j)$ the lower central series of~$\mathfrak{k}$. Then

{\rm(A)}~the image of $\mathfrak{n}_j$ under $\mathfrak{g}\to \mathfrak{g}/\mathfrak{h}$ is contained in~$\mathfrak{k}_j$;

{\rm(B)}~$\mathfrak{n}_j\to \mathfrak{k}_j$ induces a homomorphism
$\mathfrak{n}_j/(\mathfrak{h}\cap\mathfrak{n}_j)\cong \mathfrak{k}_j$.
\end{lemma}

\begin{proof}
(A)~Since $\mathfrak{h}$ is soluble, we have $\mathfrak{h}\subset\mathfrak{r}$. Since $\mathfrak{n}_1=\mathfrak{n}=[\mathfrak{g}, \mathfrak{r}]$, $\mathfrak{k}_1=\mathfrak{k}=[\mathfrak{g}/\mathfrak{h}, \operatorname{rad} (\mathfrak{g}/\mathfrak{h})]$, and $\operatorname{rad} (\mathfrak{g}/\mathfrak{h})$ is isomorphic to $\mathfrak{r}/\mathfrak{h}$ by Part~(A) of Lemma~\ref{l4.10}, it follows that the image of $\mathfrak{n}_1$ is contained in~$\mathfrak{k}_1$. The proof can be easily completed by induction.

(B)~The second isomorphism theorem implies that $\mathfrak{n}_j/(\mathfrak{h}\cap\mathfrak{n}_j) \cong (\mathfrak{h}+\mathfrak{n}_j)/\mathfrak{h}$. By Part~(A),  $(\mathfrak{h}+\mathfrak{n}_j)/\mathfrak{h}\subset \mathfrak{k}_j$. On the other hand, it is easy to see that $\mathfrak{k}_j\subset (\mathfrak{h}+\mathfrak{n}_j)/\mathfrak{h}$.
\end{proof}

\begin{proof}[Proof of Proposition~\ref{p4.9}]
We proceed by induction on the linear dimension of~$\mathfrak{g}$. When $\dim\mathfrak{g}=1$,  the assertion holds because  $\mathfrak{n}=0$ and any claim on an empty basis is true.

Assume now that $\dim\mathfrak{g}=m$ and the assertion holds for all triangular Lie algebras of linear dimension less than~$m$.
 Put $\mathfrak{z}:=\mathfrak{z}(\mathfrak{g})$ for brevity and consider three mutually exclusive cases:

(1)~$\mathfrak{z}\cap \mathfrak{n}\ne 0$;

(2)~$\mathfrak{z}=0$;

(3)~$\mathfrak{z}\cap \mathfrak{n}=0$ and $\mathfrak{z} \ne 0$.

First we show that in the cases~(1) and~(2)  there are  $x\in\mathfrak{n}$ and a finite-dimensional representation  $\pi$ of~$\mathfrak{g}$ such that (A1)--(A3) hold and~$\mathbb{K} x$ is an ideal of~$\mathfrak{g}$.

(1)~Suppose that $\mathfrak{z}\cap \mathfrak{n}\ne 0$. Then there is~$p\ge 2$ such that $\mathfrak{z}\cap \mathfrak{n}_{p-1}\ne\nobreak 0$ but $\mathfrak{z}\cap\mathfrak{n}_{p}= 0$. Since $\mathfrak{n}_{p}$ is an ideal of~$\mathfrak{g}$, we can consider the quotient algebra
$\mathfrak{g}/\mathfrak{n}_{p}$. Denote by $x$ an arbitrary non-zero element of $\mathfrak{z}\cap \mathfrak{n}_{p-1}$ and by $x'$ its pre-image under the projection $\mathfrak{g}\to \mathfrak{g}/\mathfrak{n}_{p}$. It is clear that $x'\ne 0$.

Put $\mathfrak{k}:=[\mathfrak{g}/\mathfrak{n}_{p},\mathfrak{g}/\mathfrak{n}_{p}]$ and consider the lower central series  $(\mathfrak{k}_j)$ of~$\mathfrak{k}$. Part~(A) of Lemma~\ref{l4.13} implies that $x'\in\mathfrak{k}_{p-1}$ and then $x'\in\mathfrak{z}(\mathfrak{g}/\mathfrak{n}_{p})$. Therefore $\mathfrak{z}(\mathfrak{g}/\mathfrak{n}_{p})\cap\mathfrak{k}_{p-1}\ne 0$. Moreover, it follows from Part~(B) of Lemma~\ref{l4.13} that $\mathfrak{k}_{p}=0$. So  applying Lemma~\ref{l4.11} we obtain a finite-dimensional representation of $\mathfrak{g}/\mathfrak{n}_{p}$ that satisfies (A1)--(A3) together with~$x'$. By Lemma~\ref{l4.12}, this representation can be lifted to a representation $\pi$ of $\mathfrak{g}$ that  together with~$x$ satisfies the same conditions. Also, $\mathbb{K} x$ is an ideal because $x\in\mathfrak{z}$.

(2)~Suppose that  $\mathfrak{z}=0$. It is sufficient to consider  the case when $\mathfrak{n}\ne 0$. Let $p$ denote the nilpotency degree of $\mathfrak{n}$. Then $\mathfrak{n}_{p-1}$ is a non-trivial ideal of~$\mathfrak{g}$, that is, $\mathfrak{n}_{p-1}$ is an invariant subspace of the adjoint representation~$\operatorname{ad} $ of~$\mathfrak{g}$. Denote the corresponding subrepresentation of~$\mathfrak{g}$ by~$\operatorname{ad} |_{\mathfrak{n}_{p-1}}$.

By the hypotheses, $\mathfrak{g}$ is a triangular Lie algebra. In particular, for each $y\in \mathfrak{g}$ all the eigenvalues of the linear operator~$\operatorname{ad} y$ belong to~$\mathbb{K}$ and this also obviously holds for its restriction $\operatorname{ad } |_{\mathfrak{n}_{p-1}} y$. By a generalization of Lie's theorem \cite[\S\,1.2, p.\,11, Theorem~1.2]{31}, there is a common eigenvector of  all $\operatorname{ad}|_{\mathfrak{n}_{p-1}}y$, where $y\in \mathfrak{g}$. This means that there is a non-zero $x\in \mathfrak{n}_{p-1}$ such that for every $y\in \mathfrak{g}$ there exists $\mu_y\in\mathbb{K}$ satisfying  $[y,x]=\mu_y x$. It follows that $\mathbb{K} x$ is an ideal.

We claim that $x$ and $\operatorname{ad} $ satisfy (A1)--(A3). Indeed, $\operatorname{ad} $ is injective because  $\mathfrak{z}=0$. In particular, $\operatorname{ad} x\ne 0$, that is, (A1) holds. Take $y\in\mathfrak{g}$. Since $\operatorname{ad} x(z)=-\mu_z x$ for every $z\in\mathfrak{g}$, we have
$$
\operatorname{ad} y(\operatorname{ad} x(z))=[y,[x,z]]=-\mu_z[y,x]=-\mu_y\mu_z x=\mu_y \operatorname{ad} x(z).
$$
Since~$z$ is arbitrary, $(\operatorname{ad} y)(\operatorname{ad} x)=\mu_y \operatorname{ad} x$, that is, (A2) holds. Moreover, if $y\in\mathfrak{n}$, then $[y,x]=0$ because $x\in \mathfrak{n}_{p-1}$. Hence $(\operatorname{ad} y)(\operatorname{ad} x)\,{=}\,0$, that is,
(A3) holds. Thus $x$ and $\operatorname{ad} $ satisfy (A1)--(A3).

(1)+(2) Further, we use the same argument in the cases~(1) and~(2). Denote the ideal $\mathbb{K} x$ by $\mathfrak{h}$. Lemma~\ref{l4.10} implies that $\mathfrak{g}/\mathfrak{h}$  is also triangular. It follows from Part~(B) of Lemma~\ref{l4.13} that $[\mathfrak{g}/\mathfrak{h},\mathfrak{g}/\mathfrak{h}]\cong \mathfrak{n}/\mathfrak{h}$. It is clear that $\dim\mathfrak{g}/\mathfrak{h}=m-1$ and $\dim[\mathfrak{g}/\mathfrak{h},\mathfrak{g}/\mathfrak{h}]=m-k-1$. By the induction hypothesis, there is a basis $e'_{k+1},\dots, e'_{m-1}$ of~$[\mathfrak{g}/\mathfrak{h},\mathfrak{g}/\mathfrak{h}]$ and representations $\pi'_{k+1},\dots, \pi'_{m-1}$ of $\mathfrak{g}/\mathfrak{h}$ such that (A1)--(A3) hold and  $\pi'_r(e'_j)=0$ when $r<j\le m-1$. Consider the corresponding pre-images $e_{k+1},\dots, e_{m-1}$ of these vectors in~$\mathfrak{g}$. Since $x\in \mathfrak{n}$, they are also in  $\mathfrak{n}$. Therefore $x$ complements $e_{k+1},\dots, e_{m-1}$  to a basis of~$\mathfrak{n}$. Put $\pi_m\!:=\pi$ in the case~(1) and $\pi_m\!:=\operatorname{ad} $ in the case~(2), and also $e_m\!:=x$ in both cases. Then  Lemma~\ref{l4.12} implies that~$e_j$ and the representation $\pi_j:=\pi'_j\sigma$ (where $\sigma\colon\mathfrak{g}\to\mathfrak{g}/\mathfrak{h}$) satisfy  (A1)--(A3) for $j=k+1,\dots,m-1$. Moreover, $e_m$ and $\pi_m$ satisfy (A1)--(A3) as shown above. It follows from the induction hypothesis that $\pi_r(e_j)=0$ when $r<j\le m-1$ and we have, by construction, that $\pi_r(e_m)=0$ for all $r<m$. Thus the proof of the proposition is complete in the cases~(1) and~(2).

(3)~Suppose that $\mathfrak{z}\cap \mathfrak{n}=0$ and $\mathfrak{z} \ne 0$. Then $\dim\mathfrak{g}/\mathfrak{z}<m$. It follows from Lemma~\ref{l4.10} that $\mathfrak{g}/\mathfrak{z}$ is triangular. By the induction hypothesis, there is  a basis of the nilpotent radical of~$\mathfrak{g}/\mathfrak{z}$ and representations of~$\mathfrak{g}/\mathfrak{z}$ satisfying (A1)--(A3) and the additional condition. On the one hand, the nilpotent radical of~$\mathfrak{g}/\mathfrak{z}$ coincides with $[\mathfrak{g}/\mathfrak{z},\mathfrak{g}/\mathfrak{z}]$.  On the other hand, putting $\mathfrak{h}=\mathfrak{z}$ and $j=1$ in~Part~(B) of Lemma~\ref{l4.13}, we have an isomorphism  $[\mathfrak{g}/\mathfrak{z},\mathfrak{g}/\mathfrak{z}]\cong \mathfrak{n}$. (In particular, representations just constructed can be numbered by the integers from $k+1$ to~$m$.) It follows from Lemma~\ref{l4.12}  that, after passing to pre-images, we obtain a basis of~$\mathfrak{n}$ and a tuple  $\pi_{k+1},\dots,\pi_m$ of representations of~$\mathfrak{g}$ satisfying (A1)--(A3). It is easy to see that the additional condition $\pi_r(e_j)=0$ also holds for $r<j$.
\end{proof}

In the particular case of nilpotent Lie algebras, we can show more. Namely, the following result holds. We do not need it  to prove Theorem~\ref{t4.3}  but it is useful in~\S\,\ref{ss5.1}.

\begin{proposition}
\label{p4.14}
Let~$\mathfrak{g}$ be a nilpotent Lie algebra over~$\mathbb{K}$.
Then there are a linear basis $e_{k+1},\dots, e_m$ of~$\mathfrak{n}$ and a tuple $\pi_{k+1},\dots,\pi_m$  of finite-dimensional representations of~$\mathfrak{g}$ such that $e_r$ and~$\pi_r$ satisfy \emph{(A1)}--\emph{(A3)} for every  $r\in \{k+1,\dots, m\}$ and, moreover, $\pi_r(e_j)=0$ when $r<j$. Furthermore, $\pi_r(\mathfrak{g})$ consists of nilpotent operators for each~$r$.
\end{proposition}

\begin{proof}
We use the same reasoning as in Proposition~\ref{p4.9} with the following refinements. Recall that in the proof of Proposition~\ref{p4.9} we consider  three cases (1)--(3). In the first case we refer to Lemma~\ref{l4.11}, in the proof of which it has already been established that all $\mu_y=0$ and all the operators of the representation are nilpotent. In the case~(2) we use the adjoint representation $\operatorname{ad} $. Since $\mathfrak{g}$ is nilpotent, so are all  operators $\operatorname{ad} y$. This implies that all eigenvalues are equal to~$0$. In particular, we always have that $\mu_y=0$. In the~case~(3), as well as in those places of the argument in the cases~(1) and~(2), where we  lift representations from a quotient algebra, the values of the coefficients are not changed, which means that they are also equal to~$0$. The nilpotency of the representation is also preserved.
\end{proof}

Let $\mathfrak{g}$ be a triangular Lie algebra over~$\mathbb{K}$. Fix a basis  $e_{k+1},\dots,e_m$ of~$\mathfrak{n}$ and representations $\pi_{k+1},\dots,\pi_m$ of $\mathfrak{g}$ provided by Proposition~\ref{p4.9} and complement $e_{k+1},\dots,e_m$ to a basis of~$\mathfrak{g}$ by elements $e_1,\dots, e_k$. For a multi-index $\beta=(\beta_{k+1},\dots,\beta_m)$ in $\mathbb{Z}_+^{m-k}$ denote by $\pi_\beta$ the representation of Lie algebra~$\mathfrak{g}$ (and the corresponding representation of $U(\mathfrak{g})$) that is defined as the tensor products of $\pi_{k+1},\dots,\pi_m$ taken with the multiplicities  $\beta_{k+1},\dots,\beta_m$, respectively. In particular, if $l$ and $t$ are the minimum and the maximum of  $r$ such that $\beta_r>0$, then
\begin{equation}
\label{eq4.22}
\pi_\beta(e_j):=\pi_{l}(e_j)\otimes 1\otimes\dots\otimes 1+\dots+ 1\otimes\dots\otimes 1 \otimes \pi_t(e_j),
\end{equation}
where for every~$r$ the sum has $\beta_r$ terms with the tensor factor $\pi_r(e_j)$. Moreover, $\pi_r(e_j)=0$ when $r<j$. Therefore, when $j>l$, the first terms  vanish  and we have
\begin{equation}
\label{eq4.23}
\pi_\beta(e_j)=1\otimes\dots \otimes1\otimes\pi_{i}(e_j)\otimes 1\otimes\dots\otimes 1+\dots+ 1\otimes\dots\otimes 1 \otimes \pi_t(e_j),
\end{equation}
where $i$ is the minimum number such that $i\ge j$ and  $\beta_i>0$. If $\beta=0$, we take the trivial representation as~$\pi_\beta$, that is, $\pi_\beta(e_j)=0$.

In the following two lemmas we use Condition~(A3).

\begin{lemma}
\label{l4.15}
Let $\beta=(\beta_{k+1},\dots,\beta_m)\in\mathbb{Z}_+^{m-k}\setminus\{0\}$, $l$ and $t$ be the minimum and the maximum of $r$ such that $\beta_r>0$, and $\beta_j>0$  for a certain $j\in \{l,\dots,t\}$.
Then
$$
\pi_\beta(e_j^{\beta_j}\cdots e_t^{\beta_t})=\beta!\,( 1\otimes\pi_j(e_j)\otimes \dots\otimes \pi_t(e_t)),
$$
where the tensor factors repeat with the multiplicities $\beta_j,\dots,\beta_t$ and $1$ on the left denotes for brevity the tensor product of identities.
\end{lemma}

\begin{proof}
We proceed by induction in reverse order from~$t$ to~$j$. Since $\pi_i$ and $e_i$ are taken as in Proposition~\ref{p4.9}, we have $\pi_r(e_t)=0$ when $r<t$. It follows from  $\beta_t>0$ that by~\eqref{eq4.23} $\pi_\beta(e_t)$ equals the sum
$$
1\otimes\dots\otimes \pi_t(e_t)\otimes\dots \otimes 1 +\cdots+1\otimes\dots\otimes 1 \otimes \pi_t(e_t)
$$
that have $\beta_t$ terms. Since $\pi_t(e_t)^2=0$ by Condition~(A3),  we can apply the multinomial formula and get that
$$
\pi_\beta(e_t^{\beta_t})=\beta_t!\,( 1\otimes\pi_t(e_t)\otimes \dots\otimes \pi_t(e_t)),
$$
that is,
the assertion of the lemma holds when~$j=t$.

Now assume that $s\in \{l,\dots,t-1\}$, $\beta_s>0$, and the assertion holds for all numbers greater than~$s$. In particular,
\begin{equation}
\label{eq4.24}
\pi_\beta(e_{s'}^{\beta_{s'}}\cdots e_t^{\beta_t})=\beta_{s'}!\cdots \beta_t!\,\pi_\beta(e_{s'})^{\beta_{s'}}(1\otimes\dots \otimes\pi_{s'}(e_{s'})\otimes \dots\otimes \pi_t(e_t)),
\end{equation}
where is $s'$ the least number such that $s'>s$ and $\beta_{s'}>0$.

It follows from~\eqref{eq4.23} that $\pi_\beta(e_s)=1\otimes x\otimes 1+1\otimes1\otimes y$, where
$$
x=\pi_s(e_s)\otimes 1 \otimes \dots \otimes 1 +\dots+1\otimes\dots\otimes 1 \otimes \pi_s(e_s)
$$
and $y$ contains only terms with~$\pi_r(e_s)$ when $r\ge s'$. Then by~(A3),
$$
y(\pi_{s'}(e_{s'})\otimes \dots\otimes \pi_t(e_t))=0.
$$
Hence the right-hand part in~\eqref{eq4.24} equals
$$
\beta_{s'}!\cdots \beta_t!\,(1\otimes x^{\beta_s}\otimes 1)(1\otimes\dots \otimes\pi_{s'}(e_{s'})\otimes \dots\otimes \pi_t(e_t)).
$$
Arguing as in the case when $s=t$, we conclude that
$$
x^{\beta_s}=\beta_s!\, (1\otimes\pi_s(e_s)\otimes \dots\otimes \pi_s(e_s)).
$$
Multiplying the monomials, we see that the assertion of the lemma holds for~$s$, which completes the induction.
\end{proof}

We order $\mathbb{Z}_+^{m-k}$ as follows. Put
\begin{equation}
\label{eq4.25}
\alpha\succeq \beta,\text{ if } \alpha=\beta\text{ or }\exists\,j \text{ such that } \alpha_j>\beta_j \text{ and } \alpha_i=\beta_i \ \forall\, i>j.
\end{equation}
If we read words over the alphabet $\mathbb{Z}_+$  from right to left, then this order coincides with the lexicographic order. It is sometimes called colexicographic.

\begin{lemma}[{\rm cf.~\eqref{eq4.18}}]
\label{l4.16}
Let $\beta=(\beta_{k+1},\dots,\beta_m)\in\mathbb{Z}_+^{m-k}\setminus\{0\}$. Then

{\rm(A)}~$\pi_\beta(e^\beta)\ne 0$;

{\rm(B)}~$\pi_\beta(e^\alpha)=0$ if $\alpha\succ \beta$.
\end{lemma}

\begin{proof}
Part~(A) follows from Lemma~\ref{l4.15}.

(B)~Let $\alpha\succ \beta$. Then there is $j$ such that $\alpha_j>\beta_j$ and $\alpha_i=\beta_i$ for every $i>j$. It suffices to show that $\pi_\beta(e_j^{\alpha_j}\cdots e_m^{\alpha_m})=0$.

Suppose that $\beta_j>0$. Then $\pi_\beta(e_j^{\alpha_j}\cdots e_m^{\alpha_m})=\pi_\beta(e_j)^{\alpha_j-\beta_j}\pi_\beta(e_j^{\beta_j}\cdots e_m^{\beta_m})$. By Lemma~\ref{l4.15},
$$
\pi_\beta(e_j^{\beta_j}\cdots e_m^{\beta_m})=\beta!\,( 1\otimes u\otimes \cdots\otimes \pi_t(e_t)),
$$
where $u=\pi_{j}(e_{j})\otimes \cdots\otimes \pi_j(e_j)$.

Arguing in the same way as in the proof of Lemma~\ref{l4.15}, we obtain
$$
\pi_\beta(e_j)^{\alpha_j-\beta_j}\pi_\beta(e_j^{\beta_j}\cdots e_m^{\beta_m})=\beta!\,(1\otimes x^{\alpha_j-\beta_j}\otimes 1)(1\otimes u\otimes \cdots\otimes \pi_t(e_t)),
$$
where
$$
x=\pi_j(e_j)\otimes1\otimes \dots \otimes 1 +\cdots+1\otimes\dots\otimes 1 \otimes \pi_j(e_j).
$$
Since $\pi_j(e_j)^2=0$, we have $xu=0$. Therefore  $\pi_\beta(e_j^{\alpha_j}\cdots e_m^{\alpha_m})=0$.

When $\beta_j=0$, the argument is similar with the difference that $\pi_j$ should be replaced by $\pi_{j'}$ in~$u$ and $x$, where $j'>j$ and $\beta_{j'}>0$. (If such $j'$ does not exist, then $\pi_\beta(e_{j+1}^{\beta_j}\cdots e_m^{\beta_m})=\pi_\beta(1)=0$.) Since $\pi_{j'}(e_j)\pi_{j'}(e_{j'})=0$ by Condition~(A3), we conclude that $xu=0$ and so $\pi_\beta(e_j^{\alpha_j}\cdots e_m^{\alpha_m})=0$.
\end{proof}

For every $\beta\in\mathbb{Z}_+^{m-k}$ consider the homomorphism
$$
\widetilde\pi_\beta\colon U(\mathfrak{g})\to\mathbb{K}[\lambda_1,\dots,\lambda_k]\otimes \operatorname{End} L_\beta,
$$
defined as in~\eqref{eq4.5} (here $L_\beta$ is the space of the representation $\pi_\beta$).
For $f=\sum_\alpha c_\alpha \lambda_1^{\alpha_1}\cdots\lambda_k^{\alpha_k} \in \mathbb{K}[\lambda_1,\dots,\lambda_k]$ and $\mu=(\mu_1,\dots,\mu_k)\in\mathbb{K}^k$ we define a shifted polynomial by
$$
S_\mu f:=\sum_\alpha c_\alpha(\lambda_1+\mu_1)^{\alpha_1}\cdots (\lambda_k+\mu_k)^{\alpha_k}.
$$

In the following lemma we use Condition~(A2).

\begin{lemma}[{\rm cf. Example~\ref{e4.7}}]
\label{l4.17}
Let $\beta\in \mathbb{Z}_+^{m-k}$, $\pi_\beta$ be a finite-dimensional representation of $U(\mathfrak{g})$ defined as in~\eqref{eq4.22}, $\widetilde\pi_\beta$ be the corresponding homomorphism defined as in~\eqref{eq4.5}, and $\Phi\colon \mathbb{K}[\lambda_1,\dots,\lambda_k]\to U(\mathfrak{g})$ be the ordered functional calculus corresponding to the elements $e_1,\dots,e_k$ of the basis. Then there is $\mu\in\mathbb{K}^k$ such that
$$
\widetilde\pi_\beta(\Phi(f)e^\beta)=S_\mu f\otimes \pi_\beta(e^\beta)
$$
for every $f\in \mathbb{K}[\lambda_1,\dots,\lambda_k]$.
\end{lemma}

\begin{proof}
If $\beta=0$, the assertion holds with~$\mu\,{=}\,0$. 

Let $\beta\ne0$. Suppose that $f$ is a monomial, that is, $f(\lambda_1,\dots,\lambda_k)=\lambda_1^{\gamma_1}\cdots \lambda_k^{\gamma_k}$ for some $\gamma_1,\dots,\gamma_k\in\mathbb{Z}_+$. It follows from~\eqref{eq4.5} that
\begin{align*}
\widetilde\pi_\beta(\Phi(f)e_{\beta}) &=\widetilde\pi_\beta(e_1^{\gamma_1}\cdots e_k^{\gamma_k}e_{k+1}^{\beta_{k+1}}\cdots e_m^{\beta_m})
\\
&=\prod_{j=1}^k(\lambda_j\otimes 1+1\otimes \pi_\beta(e_j))^{\gamma_j}\prod_{j=k+1}^m(1\otimes \pi_\beta(e_j)^{\beta_j}).
\end{align*}

Since Condition~(A2) is satisfied, it follows that when $j=1,\dots,k$ and $i=k+1,\dots,m$ there is $\mu_{ji} \in \mathbb{K}$ such that $\pi_i(e_j)\pi_i(e_i)=\mu_{ji} \pi_i(e_i)$. Let $l$ and $t$ be the minimum and the maximum of $r$ such that $\beta_r>0$. Since by Lemma~\ref{l4.15},
$$
\pi_\beta(e^\beta)=\beta!\,(1\otimes \pi_l(e_l)\otimes \cdots\otimes \pi_t(e_t)),
$$
where the factors repeated with the multiplicities $\beta_l,\dots,\beta_t$ and $\pi_\beta(e_j)$ is defined in~\eqref{eq4.22}, we have for every~$j$ that
$$
\pi_\beta(e_j)\pi_\beta(e^\beta)=\mu_j\pi_\beta(e^\beta),
$$
where $\mu_j:=\sum_{\beta_j>0}\mu_{ji}$. Therefore,
$$
\widetilde\pi_\beta(\Phi(f)e_{\beta})=\prod_{j=1}^k(\lambda_j +\mu_j)^{\gamma_j} \otimes \pi_\beta(e^\beta)=S_\mu f \otimes \pi_\beta(e^\beta),
$$
where $\mu=(\mu_1,\dots,\mu_k)$. The general case follows from the linearity of our maps.
\end{proof}

\subsection{The end of the proof of Theorem~\ref{t4.3}} 
\label{ss4.6}
Fix a linear basis $e_{k+1},\dots, e_m$ of~$\mathfrak{n}$ and a tuple of finite-dimensional representations  $\pi_{k+1},\dots,\pi_m$ given by Proposition~\ref{p4.9}. Complement the basis in an arbitrary way to a basis of~$\mathfrak{g}$ and consider the family  $(\pi_\beta)$ of finite-dimensional representations  of $U(\mathfrak{g})$, which  consists of tensor product and is defined in~\eqref{eq4.22}.

We continue the argument started in \S\,\ref{ss4.3}. It remains to prove that the homomorphism $\rho$ defined in~\eqref{eq4.6} is topologically injective.

The roles that played by the elements $(e_j)$ in our reasoning  are different is the cases when~$j\le k$ and when $j>k$. Write every  $a\in U(\mathfrak{g})$  as an expansion in the  elements of the PBW-basis, that is,
$$
a=\sum_{\alpha\in\mathbb{Z}_+^{m-k}} \Phi(f_\alpha)e^\alpha,
$$
where $f_\alpha\in \mathbb{R}[\lambda_1,\dots,\lambda_k]$ and $\Phi\colon \mathbb{R}[\lambda_1,\dots,\lambda_k]\to U(\mathfrak{g})$ is the ordered calculus defined in~\eqref{eq1.1}. It is obvious that this expansion is unique.

Recall that the topology on $U(\mathfrak{g})$ inherited from  $C^\infty_\mathfrak{g}$ is determined by the family
\begin{equation}
\label{eq4.26}
|a|_{\beta,M,l}:= |f_\beta|_{M,l}\qquad (l\in\mathbb{Z}_+,\quad \beta \in \mathbb{Z}_+^{m-k},\quad M\subset \mathbb{R}^m\text{ and is compact}),
\end{equation}
where  $|\,{\cdot}\,|_{M,l}:=\|\,{\cdot}\,\|_{1,M,l}$ as above. Define seminorms $\|\,{\cdot}\,\|_{\beta,M,l}$ in the same way as in~\eqref{eq2.2} using a fixed norm $\|\,{\cdot}\,\|_\beta$ on $\mathrm{T}_{d_\beta}$ instead of $\|\,{\cdot}\,\|_p$. To show that~$\rho$ is topologically injective it suffices to show that for every  $\beta$, $M$ and $l$ the seminorm $|\,{\cdot}\,|_{\beta,M,l}$ is majorized by a family $(\|\widetilde\pi_\alpha(\,{\cdot}\,)\|_{\alpha,K,n})$, where $n\in\mathbb{Z}_+$, $\alpha \in \mathbb{Z}_+^{m-k}$,
$K\subset \mathbb{R}^m$ and is compact.

Consider the order on $\mathbb{Z}_+^{m-k}$ defined in~\eqref{eq4.25}. Then $\mathbb{Z}_+^{m-k}$ is well ordered, in particular, well founded (every non-empty subset has a minimum element) and so we can proceed by induction.

Let $\beta=0$. Then $\pi_\beta$ is the trivial representation. In this case,
$$
|a|_{0,M,l}=|f_0|_{M,l}=\|\widetilde\pi_0(a)\|_{0,M,l},
$$
and the assertion holds.

Now let $\beta\in\mathbb{Z}_+^{m-k}\setminus\{0\}$. Suppose that for every $\gamma\prec\beta$, $M$ and $l$ the seminorm $|\,{\cdot}\,|_{\gamma,M,l}$ is majorized by $(\|\widetilde\pi_\alpha(\,{\cdot}\,)\|_{\alpha,K,n})$. Further, fix $M$ and $l$. By Lemma~\ref{l4.17}, there is a vector~$\mu$ in~$\mathbb{R}^k$ such that
\begin{equation}
\label{eq4.27}
\widetilde\pi_\beta(\Phi(f)e^\beta)=S_\mu f\otimes \pi_\beta(e^\beta)
\end{equation}
for all~$f\in \mathbb{K}[\lambda_1,\dots,\lambda_k]$. Denote the set
$$
\{(r_1-\mu_1,\dots,r_k-\mu_k)\colon (r_1,\dots,r_k)\in M\}
$$
by $M-\mu$ and put $C\!:=\|\pi_\beta(e^\beta)\|_\beta$. (Recall that $\|\,{\cdot}\,\|_\beta$ is a fixed submultiplicative norm on $\mathrm{T}_{d_\beta}$, where $d_\beta$ is the dimension of~$\pi_\beta$.)

It is easy to see that $|a|_{\beta,M,l}=|f_\beta|_{M,l}=|S_\mu f_\beta|_{M-\mu,l}$. On the other hand, it follows from~\eqref{eq4.27} and the definition of $\|\,{\cdot}\,\|_{\beta,M-\mu,l}$ that
$$
\|\widetilde\pi_\beta(\Phi(f)e^\beta)\|_{\beta,M-\mu,l}=C\, |S_\mu f|_{M-\mu,l}
$$
for every~$f$. Lemma~\ref{l4.16} implies that $\pi_\beta(e^\beta)\ne 0$. So since $\|\,{\cdot}\,\|_\beta$ is a norm, we have that $C\ne 0$. Thus,
\begin{equation}
\label{eq4.28}
|a|_{\beta,M,l}=C^{-1}\,\|\widetilde\pi_\beta(\Phi(f_\beta)e^\beta)\|_{\beta,M-\mu,l}.
\end{equation}

Consider the projection on $C^\infty_\mathfrak{g}$ given by
$$
P_\beta\biggl(\sum_{\alpha} \Phi(f_\alpha)e^\alpha\biggr):=\sum_{\alpha\prec \beta} \Phi(f_\alpha)e^\alpha.
$$
It follows from Part~(B) of Lemma~\ref{l4.16} that $\pi_\beta(e^\alpha)=0$ when $\alpha\succ \beta$. Since the order on $\mathbb{Z}_+^{m-k}$ defined above is linear, we have that
$$
\widetilde\pi_\beta(\Phi(f_\beta)e^\beta)=\widetilde\pi_\beta(a)-\widetilde\pi_\beta P_\beta(a)
$$
for every $a\in C^\infty_\mathfrak{g}$. Then~\eqref{eq4.28} implies that
\begin{equation}
\label{eq4.29}
|a|_{\beta,M,l}\le C^{-1}(\|\widetilde\pi_\beta(a)\|_{\beta,M-\mu,l}+\|\widetilde\pi_\beta P_\beta(a)\|_{\beta,M-\mu,l})
\end{equation}
(cf.~\eqref{eq4.12} and~\eqref{eq4.17}). Note that $\widetilde\pi_\beta$ is continuous and so $\|\widetilde\pi_\beta P_\beta(\,{\cdot}\,)\|_{\beta,M-\mu,l}$ is majorized by the family $(|P_\beta(\,{\cdot}\,)|_{\beta',M',l'})$. Since $|P_\beta(a)|_{\beta',M',l'}=|a|_{\beta',M',l'}$ when $\beta'\prec\beta$ and $|P_\beta(a)|_{\beta',M',l'}=0$ when $\beta'\succeq \beta$ for every $a\in C^\infty_\mathfrak{g}$, we have by the induction hypothesis that $\|\widetilde\pi_\beta P_\beta(\,{\cdot}\,)\|_{\beta,M-\mu,l}$ is majorized by $(\|\widetilde\pi_\alpha(\,{\cdot}\,)\|_{\alpha,K,n})$. It follows that  $|\,{\cdot}\,|_{\beta,M,l}$ is majorized by the last family.

Thus the principle of induction implies that $|\,{\cdot}\,|_{\beta,M,l}$ is majorized by  $(\|\widetilde\pi_\alpha(\,{\cdot}\,)\|_{\alpha,K,n})$.
So~$\rho$ is topologically injective, which means that the assertion of Theorem~\ref{t4.3} holds in the case of the~special basis of~$\mathfrak{n}$ provided by Proposition~\ref{p4.9}.  To complete the proof of the theorem, we recall that the independence of the topology on $C^\infty_\mathfrak{g}$ from the choice of a basis of~$\mathfrak{n}$ is shown in the first part of the argument (see \S\,\ref{ss4.3}). The proof of Theorem~\ref{t4.3} is complete.

\subsection{Proof of the theorem on multiplicative $C^\infty_\mathfrak{g}$-functional calculus}
\label{ss4.7}

The argument invokes the main results of this section and \S\,\ref{s3}.

\begin{proof}[Proof of Theorem~\ref{t4.4}]
Let $\mathfrak{g}$ be a triangular finite-dimensional real Lie algebra, $B$ a projective limit of real Banach algebras of polynomial growth and $\gamma\colon \mathfrak{g}\to B$ a homomorphism of real Lie algebras. It suffices to consider the case when $B$ is a Banach algebra.

Note that $\gamma$ and the embedding $\mu\colon \mathfrak{g}\to C^\infty_\mathfrak{g}$ extend to homomorphisms with domain $U(\mathfrak{g})$. Let $e_{k+1},\dots, e_m$ be a linear basis of~$[\mathfrak{g},\mathfrak{g}]$ and $e_1,\dots, e_k$ a complement to a linear basis of~$\mathfrak{g}$. Then $\gamma(e_{k+1}),\dots,\gamma(e_m)\in [B,B]$. By Theorem~\ref{t2.8}, these elements are in $\operatorname{Rad} B$ and so they are topologically nilpotent. Moreover, having polynomial growth, they are nilpotent by Proposition~\ref{p2.4}. Since $\gamma(e_1),\dots,\gamma(e_k)$ are of polynomial growth, it follows from Theorem~\ref{t3.3} that there is a continuous linear map $\theta\colon C^\infty_\mathfrak{g}\to B$ such that $\gamma=\theta\mu$ (recall that, by definition, $C^\infty_\mathfrak{g}=C^\infty(\mathbb{R}^k)\,{\mathbin{\widehat{\otimes}}}\, \mathbb{R}[[e_{k+1},\dots,e_m]]$ 
as a locally convex space). By Theorem~\ref{t4.3}, the multiplication in~$C^\infty_\mathfrak{g}$ is a continuous extension of the multiplication in~$U(\mathfrak{g})$ and so  $\theta$ is a homomorphism. Thus the universal property is proved. It immediately follows that the algebra $C^\infty_\mathfrak{g}$ is independent of the choice of a basis of~$\mathfrak{n}$ and a complement to a basis of~$\mathfrak{g}$.
\end{proof}

\begin{remark}
\label{r4.18}
Note that the hypotheses of Theorem~\ref{t4.4} can be weakened by requiring only that the elements $\gamma(e_1),\dots,\gamma(e_m)$ are of polynomial growth (without assuming this for all elements of $B$). Under such weakened hypotheses, we can nevertheless claim that $\theta$ is well defined, and, moreover, $\theta(f)$ is of polynomial growth for any $f\in C^\infty_\mathfrak{g}$. Indeed, in the proof of Theorem~\ref{t4.4} we use only the facts that  $\gamma(e_1),\dots,\gamma(e_m)$ are of polynomial growth and $\gamma(e_{k+1}),\dots,\gamma(e_m)$ are topologically nilpotent. To check the latter condition, extend $\gamma$ to a homomorphism $\mathfrak{g}_\mathbb{C}\to B_\mathbb{C}$ between the complexifications. Since $\mathfrak{g}_\mathbb{C}$ is soluble, it follows from results of Turovskii \cite{35} (see also a proof in~\cite[\S\,24, p.\,130, Theorem~1]{29}) that $[\gamma(\mathfrak{g}_\mathbb{C}),\gamma(\mathfrak{g}_\mathbb{C})]$ is contained in the radical of $B_0$, the closed subalgebra of~$B_\mathbb{C}$ generated by $\gamma(\mathfrak{g}_\mathbb{C})$. This implies that the elements of interest are topologically nilpotent. Moreover, by Theorem~\ref{t4.3}, every $f\in C^\infty_\mathfrak{g}$ is of polynomial growth. Since the property of having polynomial growth is preserved under continuous homomorphisms of Arens--Michael algebras, $\theta(f)$ also is of polynomial growth.
\end{remark}

The question of whether the hypotheses of Theorem~\ref{t4.4} can be further weakened by assuming polynomial growth only for the images of the algebraic generators of~$\mathfrak{g}$ (that is, for $\gamma(e_1) ,\dots,\gamma(e_k)$) remains open. It can be reformulated as follows.

\begin{question}
\label{q4.19}
Assume that a finite-dimensional real Lie subalgebra~$\mathfrak{h}$ of a~real Banach algebra is generated by a finite set of elements of polynomial growth. Does it follow that all elements~$\mathfrak{h}$ are of polynomial growth?
\end{question}

\section{Locally defined non-commutative functions of class $C^\infty$ and sheaves}
\label{s5}

\subsection{Nilpotent case}
\label{ss5.1}
The theorems in this section are analogues of Dosi's results \cite{7}, \cite{2} on algebras of non-commutative holomorphic functions; for details see \S\,\ref{s6}.

Suppose that $\mathfrak{g}$ is nilpotent. We now define algebras of non-commutative smooth functions on open subsets of $\mathbb{R}^k$. (Here $k$ denotes the dimension of~$\mathfrak{g}/\mathfrak{n}$ and, in addition,  $m$ denotes the dimension of $\mathfrak{g}$; we keep the notation in Theorem~\ref{t4.3}.) As above, we fix a linear basis $e_{k+1},\dots, e_m$ of~$\mathfrak{n}$  and a complement $e_1,\dots, e_k$  to a linear basis of~$\mathfrak{g}$. For an open subset $V$ of $\mathbb{R}^k$
we consider the Fr\'echet space
\begin{equation}
\label{eq5.1}
C^\infty_\mathfrak{g}(V):=C^\infty(V)\mathbin{\widehat{\otimes}} \mathbb{R}[[e_{k+1},\dots,e_m]]
\end{equation}
and identify $U(\mathfrak{g})$ with a subspace of $C^\infty_\mathfrak{g}(V)$ (which is dense by the Weierstrass approximation  theorem for continuously differentiable functions \cite[Theorem 1.6.2]{36}). In particular, we put $C^\infty_\mathfrak{g}(\emptyset)=0$. It is clear that $C^\infty_\mathfrak{g}(\mathbb{R}^k)=C^\infty_\mathfrak{g}$.

\begin{theorem}
\label{t5.1}
Let $\mathfrak{g}$ be a nilpotent real Lie algebra and $V$ an open subset of $\mathbb{R}^k$. Then the multiplication in~$U(\mathfrak{g})$ extends to a continuous multiplication in $C^\infty_\mathfrak{g}(V)$ and the topology is independent on the choice of a basis of~$\mathfrak{n}$. Moreover, with this multiplication, $C^\infty_\mathfrak{g}(V)$ is a projective limit of real Banach algebras of polynomial growth and hence a Fr\'echet--Arens--Michael algebra of polynomial growth.
\end{theorem}

The difference between the proofs of Theorems~\ref{t5.1} and~\ref{t4.3} is that in the nilpotent case we estimate seminorms for fixed compact subsets of~$V$ (since the shift operations are trivial). Recall that $\|\,{\cdot}\,\|_{p,K,n}$ are defined in~\eqref{eq2.2} and $|\,{\cdot}\,|_{\beta,M,l}$ in~\eqref{eq4.26}, and each finite-dimensional representation $\pi\colon U(\mathfrak{g})\to\operatorname{End} L$ induces a homomorphism
$$
\widetilde\pi\colon U(\mathfrak{g})\to\mathbb{R}[\lambda_1,\dots,\lambda_k]\otimes \operatorname{End} L
$$
defined in~\eqref{eq4.5}.


We need an auxiliary lemma.

\begin{lemma}
\label{l5.2}
Suppose that $\mathfrak{g}$ is a nilpotent real Lie algebra. Let $n\in\mathbb{Z}_+$, $K$ be a compact subset of~$\mathbb{R}^k$, and $\pi$ be a finite-dimensional representation of $\mathfrak{g}$  whose image is a nilpotent set. Then the seminorm $\|\widetilde\pi(\,{\cdot}\,)\|_{d,K,n}$, where $d$ is the dimension of~$\pi$, is majorized by a family $(|\,{\cdot}\,|_{\beta,K,l})$, where $\beta\in \mathbb{Z}_+^{m-k}$ and $l\in\mathbb{Z}_+$.
\end{lemma}

\begin{proof}
Since $\pi(\mathfrak{g})$ is nilpotent, so is $\pi(U(\mathfrak{g}))$.
Fix $t\in\mathbb{Z}_+$ such that $\pi(U(\mathfrak{g}))^t=0$. (We could even put $t=d$ but this does not matter.) Let, as above, $\Phi\colon \mathbb{R}[\lambda_1,\dots,\lambda_k]\to U(\mathfrak{g})$ be an ordered calculus; see~\eqref{eq1.1}. Since $\lambda_i\otimes 1$ and $1\otimes \pi(e_j)$ (see~\eqref{eq4.5}) commute, we have by the Taylor formula that
$$
\widetilde\pi(\Phi(f))=\sum_{\alpha\in \mathbb{Z}_+^k,\, |\alpha|<t} \frac{f^{(\alpha)}(\lambda_1,\dots,\lambda_k)}{\alpha_1!\cdots \alpha_k!}\otimes \pi(e_1^{\alpha_1}\cdots e_k^{\alpha_k})
$$
for every $f\in \mathbb{R}[\lambda_1,\dots,\lambda_k]$. Since each element of $U(\mathfrak{g})$ has the form
$$
a=\sum_{\beta\in \mathbb{Z}_+^{m-k}} \Phi(f_\beta)e^\beta,
$$
we obtain
\begin{equation}
\label{eq5.2}
\widetilde\pi(a)=\sum_{|\beta|<t} \widetilde\pi(\Phi(f_\beta)e^\beta)=\sum_{|\alpha|+ |\beta|<t} \frac{f_\beta^{(\alpha)}(\lambda_1,\dots,\lambda_k)}{\alpha_1!\cdots \alpha_k!}\otimes \pi(e_1^{\alpha_1}\cdots e_k^{\alpha_k} e_{k+1}^{\beta_1}\cdots e_m^{\beta_{m-k}}).
\end{equation}
For every $\alpha$ and $\beta$ the value of the seminorm $\|\,{\cdot}\,\|_{d,K,n}$ at the corresponding term on the right-hand side is majorized by $|f_\beta|_{K,n+t}$. Since, by definition, $|a|_{\beta,K,n+t}= |f_\beta|_{K,n+t}$, the assertion of the lemma follows.
\end{proof}

\begin{proof}[Proof of Theorem~\ref{t5.1}]
The independence of the choice of a basis of~$\mathfrak{n}$ can be shown as in \S\,\ref{ss4.3}.

As in the proof of Theorem~\ref{t4.3}, we choose a basis $e_{k+1},\dots, e_m$ of~$\mathfrak{n}$ in a convenient way. Namely, we consider $e_{k+1},\dots, e_m$ and a tuple  $\pi_{k+1},\dots,\pi_m$ of finite-dimensional representations provided by Proposition~\ref{p4.14} and complement these vectors  to a basis of~$\mathfrak{g}$ in an arbitrary way. Consider also the family $(\pi_\beta;\, \beta\in \mathbb{Z}_+^{m-k})$ of finite-dimensional representations of $U(\mathfrak{g})$ consisting of tensor products and defined in~\eqref{eq4.22}.

Since all operators in $\pi_r(\mathfrak{g})$ are nilpotent for every~$r$, all operators in $\pi_\beta(\mathfrak{g})$ are nilpotent for every~$\beta$. Recall also that a finite-dimensional representation of a triangular (in particular, nilpotent) Lie algebra is equivalent to an upper triangular one. It follows that every operator in $\pi_\beta(\mathfrak{g})$ is a nilpotent set. Applying Lemma~\ref{l5.2}, we get that for every compact subset of~$V$ and every $\beta\in\mathbb{Z}_+^{m-k}$  the homomorphism $\widetilde\pi_\beta$ is continuous with respect to the restriction of the topology on $C^\infty_\mathfrak{g}(V)$ to~$U(\mathfrak{g})$ and the topology on $\mathbb{R}[\lambda_1,\dots,\lambda_k]\otimes \operatorname{End} L$ determined by the family $(\|\,{\cdot}\,\|_{d_\beta,K,n})$, where $d_\beta$ is the dimension of $\pi_\beta$, $K$ is a compact subset of~$V$ and $n\in\mathbb{N}$. Moreover, we can assume that $\widetilde{\pi_\beta}$ is a continuous homomorphism with values in~$C^{\infty}(V,\mathrm{T}_{d_\beta})$.

Further, we argue in the same way as in the proof of Theorem~\ref{t4.3}. Consider the following continuous homomorphism:
$$
\rho_V\colon U(\mathfrak{g})\to \prod_{\beta} C^{\infty}(V,\mathrm{T}_{d_\beta})\colon  a\mapsto (\widetilde\pi_\beta(a))
$$
(cf.~\eqref{eq4.6}). To complete the proof it suffices to show that $\rho_V$ is topologically injective. Namely, we need to check that for every $\beta\in \mathbb{Z}_+^{m-k}$, compact subset~$M$ of~$V$ and $l\in\mathbb{Z}_+$ the seminorm $|\,{\cdot}\,|_{\beta,M,l}$ is majorized by the family $(\|\widetilde\pi_\alpha(\,{\cdot}\,)\|_{\alpha,M,n})$, where $n\in\mathbb{Z}_+$ and $\alpha \in \mathbb{Z}_+^{m-k}$.

We proceed by induction in the same way as in \S\,\ref{ss4.6}. The difference is in the fact that, in the nilpotent case, Proposition~\ref{p4.14} implies that $\mu=0$ in~\eqref{eq4.27}, that is, $M$ is shifted. Therefore, the seminorm on the right-hand sides in~\eqref{eq4.28} and~\eqref{eq4.29} takes the form $\|\,{\cdot}\,\|_{\beta,M,l}$ and it follows from the corresponding estimates that $\rho_V$ is topologically injective.
\end{proof}

\begin{remark}
\label{r5.3}
According to Theorem~\ref{t5.1}, the algebra $C^\infty_\mathfrak{g}(V)$ is independent of the choice of the basis $e_{k+1},\dots, e_m$ of~$\mathfrak{n }$. However, from a~formal point of view, they may depend on the choice of the complement $e_1,\dots, e_k$ to a basis of~$\mathfrak{g}$. To show the independence of the choice of the complement one need a technique that goes beyond the scope of this article (in particular, a local universal property that generalizes Theorem~\ref{t4.4}).
\end{remark}

Next, we show that the algebras of the form $C^\infty_\mathfrak{g}(V)$ are consistent with each other, namely, they form a sheaf on the Gelfand spectrum of $C^\infty_\mathfrak{g}$, which is homeomorphic to $ \mathbb{R}^k$ (cf. \cite[\S\,5.3, p.\,123]{7} and Theorem~\ref{t6.5} on non-commutative holomorphic functions).

\begin{remark}
\label{r5.4}
Recall that a sheaf of Fr\'echet spaces (Fr\'echet algebras) is a presheaf of Fr\'echet spaces (Fr\'echet algebras\footnote{We consider all categories of locally convex algebras with continuous homomorphisms as morphisms.}) that is also a sheaf of sets (after applying the forgetful functor) \cite[\S\,6]{37}; see also \cite[\S\,4.3, p.\,109]{28} and \cite[Appendix~A]{38}. This definition is due to the fact that the gluing axiom has sense only for countable covers (a general product of Fr\'echet spaces is not necessarily a Fr\'echet space). However, it follows from the inverse mapping theorem for Fr\'echet spaces (which in this context means that the forgetful functor reflects isomorphisms; cf. the algebraic case in~\cite[Part~1, \S\,6.9]{39}) that, for countable covers, it suffices to verify the gluing axiom in the category of sets. If the base space of a sheaf admits a countable base of the topology, then it is a Lindel\"{o}f space (any open cover contains a countable subcover). In this case, it is easy to show that if the gluing axiom holds for countable covers, then it also holds for arbitrary covers but already in the category of locally convex spaces (algebras), which is wider. Thus a sheaf of Fr\'echet spaces (Fr\'echet algebras) on a space with countable base is always a sheaf of locally convex spaces (algebras) and we see that there are no problems with the above definition. Similarly, we consider sheaves of Fr\'echet--Arens-Michael algebras (of polynomial growth) defining them as presheaves of Fr\'echet--Arens--Michael algebras (of polynomial growth) that are also sheaves of sets.
\end{remark}

Let $V$ and $W$ be open subsets of~$\mathbb{R}^k$ such that $W\subset V$. Then the topology on~$U(\mathfrak{g})$ induced by the embedding $U(\mathfrak{g})\to C^\infty_\mathfrak{g}(V)$ is stronger that  topology induced by the embedding $U(\mathfrak{g})\to C^\infty_\mathfrak{g}(W)$. Since both algebras are completions of~$U(\mathfrak{g})$, we have a continuous homomorphism $\tau_{VW}\colon C^\infty_\mathfrak{g}(V)\to C^\infty_\mathfrak{g}(W)$.
\setcounter{footnote}{-1}
Recall that the Gelfand spectrum\footnote{Added in the English translation: Strictly speaking, this is the real spectrum because the Gelfand spectrum is usually defined as the set of homomorphisms to $\mathbb{C}$. But for algebras of polynomial growth these spectra coincide.}
of an Arens-Michael $\mathbb{R}$-algebra is the set of continuous homomorphisms to~$\mathbb{R}$ endowed with the  weak* topology. \addtocounter{footnote}{9}

\begin{theorem}
\label{t5.5}
Let $\mathfrak{g}$ be a nilpotent real Lie algebra. Then the correspondences
$$
V\mapsto C^\infty_\mathfrak{g}(V)\quad\text{and}\quad (W\subset V)\mapsto \tau_{VW}
$$
determine a sheaf of Fr\'echet--Arens--Michael  $\mathbb{R}$-algebras of polynomial growth on the Gel\-fand spectrum of $C^\infty_\mathfrak{g}$.
\end{theorem}

\begin{proof}
It follows from Theorem~\ref{t2.8} that the Gelfand spectrum of $C^\infty_\mathfrak{g}$ coincides with the Gelfand spectrum of  $C^\infty(\mathbb{R}^k)$, that is, $\mathbb{R}^k$. It is easy to see that we have a  contravariant functor from the category of open subsets of~$\mathbb{R}^k$ to the category of Arens--Michael  $\mathbb{R}$-algebras of polynomial growth. Thus we obtain a presheaf. To show that this presheaf is a sheaf it remains to verify the gluing axiom in the category of sets. Since for every open set $V$ the Fr\'echet space $C^\infty_\mathfrak{g}(V)$ is the product of a countable number of copies of $C^\infty(V)$ it suffices to check that the gluing axiom holds for the functions of class $C^\infty$ on~$\mathbb{R}^k$. The last fact is well known (and can be verified directly).
\end{proof}

\subsection{Discussion on the general case}
\label{ss5.2}
Let, as above, $\mathfrak{a}\mathfrak{f}_1$ denote the two-dimensional real Lie algebra with  basis $e_1,e_2$ and multiplication defined by the relation $[e_1,e_2]=e_2$. Recall that $\mathfrak{a}\mathfrak{f}_1$ is not nilpotent.

In contrast to the case of a nilpotent Lie algebra (see Theorem~\ref{t5.1}),  the multiplication in~$U(\mathfrak{a}\mathfrak{f}_1)$ cannot be extended to a continuous multiplication in~$C^\infty(V)\,{\mathbin{\widehat{\otimes}}}\, \mathbb{R}[[e_2]]$ for every open subset~$V$ of~$\mathbb{R}$.
Indeed, put $V{=}\,(0,2)$ and let  $(f_n)$ be a sequence of polynomials in~$\lambda$ converging to the function $1/\lambda$ in the space $C^\infty(0,2)$. Since $e_2f_n(e_1)=f_n(e_1-1)e_2$, the assumption  that the multiplication is continuous implies that  $f_n(e_1-1)e_2$ converges in $C^\infty(0,2)\mathbin{\widehat{\otimes}}\mathbb{R}[[e_2]]$ and so $f_n(\lambda-1)$ converges in~$C^\infty(0,2)$. But this is not the case and we get a contradiction.

Thus the Gelfand spectrum of a general triangular Lie algebra is too small to construct a sheaf of non-commutative smooth functions even in the simplest case when $\mathfrak{g}=\mathfrak{a}\mathfrak{f}_1$. However, it is possible to do it on a larger space, at least for this algebra; see \cite{40}. So, in the nilpotent case, we get sheaves that can be treated as analytic analogues of sheaves on formal non-commutative schemes considered by Kapranov in~\cite{41} but the situation is more complicated for general triangular Lie algebras.

\section{Algebras of non-commutative holomorphic functions}
\label{s6}

In~\cite{2} Dosi considered a topological algebra $\mathfrak{F}_\mathfrak{g}$ of `formally radical entire functions' associated with a basis of a positively graded nilpotent complex Lie algebra~$\mathfrak{g}$ consistent with the lower central series, and proved that the multiplication induced by that in the universal enveloping algebra $U(\mathfrak{g})$ (over the field~$\mathbb{C}$) is jointly continuous. Thus $\mathfrak{F}_\mathfrak{g}$ is a Fr\'echet algebra (since the topology can be determined by a countable family of seminorms). In what follows, arguing  as in the proof of Theorem~\ref{t4.3}, we not only show that $\mathfrak{F}_\mathfrak{g}$ is an Arens--Michael algebra but also define $\mathfrak{F}_\mathfrak{g}$ for an arbitrary soluble complex Lie algebra~$\mathfrak{g}$.

Indeed, suppose that $\mathfrak{g}$ is soluble and fix a linear basis $e_{k+1},\dots, e_m$ of~$\mathfrak{n}$ with a complement $e_1,\dots, e_k$ to a linear basis of~$\mathfrak{g}$. Consider the Fr\'echet space
$$
\mathfrak{F}_\mathfrak{g}:=\mathcal{O}(\mathbb{C}^k)\mathbin{\widehat{\otimes}} \mathbb{C}[[e_{k+1},\dots,e_m]],
$$
where $\mathcal{O}(\mathbb{C}^k)$ is the algebra of holomorphic functions on~$\mathbb{C}^k$; cf.~\eqref{eq4.1}. Then the following analogue of Theorem~\ref{t4.3} holds.

\begin{theorem}
\label{t6.1}
Let $\mathfrak{g}$ be a soluble complex Lie algebra. Then the multiplication in~$U(\mathfrak{g})$ extends to a jointly continuous multiplication in $\mathfrak{F}_\mathfrak{g}$ and the topology is independent of the choice of a basis of~$\mathfrak{n}$. Moreover, $\mathfrak{F}_\mathfrak{g}$ is a Fr\'echet--Arens--Michael $\mathbb{C}$-algebra with respect to this multiplication.
\end{theorem}

\begin{proof}
With slight modifications, we use the same argument as in Theorem~\ref{t4.3}. As in~\S\,\ref{ss4.5}, we construct a countable family  $(\pi_\beta\colon U(\mathfrak{g})\to \mathrm{T}_{d_\beta})$ of homomorphisms and the corresponding family $\widetilde\pi_\beta\colon U(\mathfrak{g})\to\mathcal{O}(\mathbb{C}^k,\mathrm{T}_{d_\beta})$ defined by~\eqref{eq4.5}. Further, we consider the homomorphism
$$
\rho\colon U(\mathfrak{g})\to \prod_{\beta} \mathcal{O}(\mathbb{C}^k,\mathrm{T}_{d_\beta})\colon  a\mapsto (\widetilde\pi_\beta(a)),
$$
similar to that in~\eqref{eq4.6}. To show that it can be extended to a continuous linear map from $\mathfrak{F}_\mathfrak{g}$ we write $\mathcal{O}(\mathbb{C}^m)$ as the projective tensor product of  $m$ copies of $\mathcal{O}(\mathbb{C})$. Applying the  theorem on holomorphic calculus to each factor, we can
extend  $\widetilde\pi_\beta$ to a continuous linear map $\mathcal{O}(\mathbb{C}^m)\to\mathcal{O}(\mathbb{C}^k,\mathrm{T}_{d_\beta})$. Since $\widetilde\pi_\beta(e_{k+1}),\dots, \widetilde\pi_\beta(e_m)$ are nilpotent, it factors through~$\mathfrak{F}_\mathfrak{g}$.

Every soluble complex Lie algebra is triangular. All the auxiliary assertions in \S\,\ref{ss4.5}  used in the proof of Theorem~\ref{t4.3} hold in the case of the field~$\mathbb{C}$. Therefore the argument in \S\,~\ref{ss4.6} works without changes for holomorphic functions. This implies that $\rho$ is topologically injective, which completes the proof.
\end{proof}

It is clear that $\mathfrak{F}_\mathfrak{g}=\mathcal{O}(\mathbb{C}^k)$ when $\mathfrak{g}$ is abelian. Thus, if $\mathfrak{g}$ is a soluble Lie algebra over~$\mathbb{C}$, then $\mathfrak{F}_\mathfrak{g}$ can be treated as an algebra of entire functions in non-commutative variables generating~$\mathfrak{g}$. It is worth noting that, generally speaking, $\mathfrak{F}_\mathfrak{g}$ does not coincide with the Arens--Michael envelope of $U(\mathfrak{g})$. The latter is described explicitly in~\cite{42} (the nilpotent case) and \cite{43} (the general case). Nevertheless,  $\mathfrak{F}_\mathfrak{g}$  satisfies a certain universal property but this topic requires separate consideration. We also note that  the concomitant  Banach algebras of $\mathfrak{F}_\mathfrak{g}$ are not of polynomial growth (since they are defined over~$\mathbb{C}$). However, such an algebra has the analogous property  of being `supernilpotent' (this means that all elements of the commutant are nilpotent) and thus satisfies a polynomial identity.

Dosi also considered a local variant of the algebra~$\mathfrak{F}_\mathfrak{g}$. Namely, let
$$
\mathfrak{F}_\mathfrak{g}(V):=\mathcal{O}(V)\mathbin{\widehat{\otimes}} \mathbb{C}[[e_{k+1},\dots,e_m]],
$$
where $\mathfrak{g}$ is a nilpotent complex Lie algebra and $V$ is an open subset of~$\mathbb{C}^k$. (The above Definition~\eqref{eq5.1} is an obvious analogue of this definition.)

Under the assumptions that $\mathfrak{g}$ is positively graded and the basis is consistent with the lower central series, Dosi showed first that for an open subset $V\subset \mathbb{C}^k$, the multiplication in~$U(\mathfrak{g})$ induces\footnote{In the case  when $V$ is a polydisc, $U(\mathfrak{g})$ is dense in~$\mathfrak{F}_\mathfrak{g}(V)$ and the multiplication can be extended by continuity. In the general case, we can use a cover by sets in a base that contains all polydiscs; cf. Theorem~\ref{t6.6}.} a continuous multiplication in $\mathfrak{F}_\mathfrak{g}(V)$ and so it is a Fr\'echet algebra (\cite[\S\,5.1, p.\,120]{7} or \cite[\S\,5.6, p.\,21]{2}) and, secondly, that the corresponding functor is a sheaf on its (complex) Gelfand spectrum \cite[\S\,5.3, p.\,123]{7}; cf. Theorem~\ref{t5.5} above.

Theorems~\ref{t6.2}, \ref{t6.5} and~\ref{t6.6} proved below are strengthenings of Dosi's results. It differs from results in~\cite{7} in three points: (1)~we do not assume that $\mathfrak{g}$ is positively graded; (2)~we do not assume that the basis is consistent with the lower central series; (3)~we prove that $\mathfrak{F}_\mathfrak{g}(V)$ is not only a Fr\'echet algebra but also an Arens--Michael algebra.

In the case when the algebra of polynomials is not dense in~$\mathcal{O}(V)$, the algebra $U(\mathfrak{g})$ is not dense in~$\mathfrak{F}_\mathfrak{g}(V)$, and we cannot extend the multiplication from~$U(\mathfrak{g})$ by  continuity. So we first consider those open subsets in~$\mathbb{C}^k$ for which polynomials are dense in~$\mathcal{O}(V)$.

\begin{theorem}
\label{t6.2}
Let $\mathfrak{g}$ be a nilpotent complex Lie algebra and $V$ an open subset in $\mathbb{C}^k$ such that the algebra of polynomials is dense in~$\mathcal{O}(V)$.
Then the multiplication in~$U(\mathfrak{g})$  extends to a continuous multiplication in $\mathfrak{F}_\mathfrak{g}(V)$. Moreover, with this multiplication, $\mathfrak{F}_\mathfrak{g}(V)$ is a Fr\'echet--Arens--Michael algebra and independent of the choice of a basis of~$\mathfrak{n}$.
\end{theorem}

\begin{remark}
\label{r6.3}
The question of the independence of algebras of the form $\mathfrak{F}_\mathfrak{g}(V)$ from the choice of the complement  $e_1,\dots, e_k$ to a basis of~$\mathfrak{g}$ remains open.
\end{remark}

To prove Theorem~\ref{t6.2} we need the following analogue of Lemma~\ref{l5.2}. Note that the topologies on $\mathcal{O}(V,\mathrm{T}_q)$ and $\mathfrak{F}_\mathfrak{g}(V)$ are determined by the families  $(\|\,{\cdot}\,\|_{q,K,0})$ and $(|\,{\cdot}\,|_{\beta,M,0})$  of seminorms, which are defined in the same way as in~\eqref{eq2.2} and~\eqref{eq4.26}, respectively, but with the difference that $K$ and $M$ are now compact subsets of $\mathbb{C}^k$ (not of $\mathbb{R}^k$) contained in~$V$.

\begin{lemma}
\label{l6.4}
Suppose that  $\mathfrak{g}$ is a nilpotent complex Lie algebra. Let $V$ be an open subset of~$\mathbb{C}^k$, $K$ a compact subset of~$V$ and $\pi$ a finite-dimensional representation of $\mathfrak{g}$ whose range is a nilpotent set. Then there is a compact set $K'$ such that $K\subset K'\subset V$ and $\|\widetilde\pi(\,{\cdot}\,)\|_{d,K,0}$, where $d$ is the dimension of~$\pi$, is majorized by the family  $(|\,{\cdot}\,|_{\beta,K',0})$, where $\beta\in \mathbb{Z}_+^{m-k}$.
\end{lemma}

\begin{proof}
The argument is similar to that for Lemma~\ref{l5.2}. The only difference is that to estimate the derivatives in~\eqref{eq5.2} we need the Weierstrass theorem on the uniform convergence of derivatives (see, for example, \cite[Chapter~1, \S\,2, Russain p.\,34, Theorem~8]{44}), which implies the existence of the desired compact subset~$K'$.
\end{proof}

\begin{proof}[Proof of Theorem~\ref{t6.2}]
We combine the arguments for Theorems~\ref{t5.1} and~\ref{t6.1}. In the same way as in the proof of Theorem~\ref{t5.1}, but  using holomorphic  functions instead of smooth (as in the proof of Theorem~\ref{t6.1}), we construct a homomorphism
$$
\rho_V\colon U(\mathfrak{g})\to \prod_{\beta} \mathcal{O}(V,\mathrm{T}_{d_\beta}).
$$
Identify $U(\mathfrak{g})$ with a dense subspace of~$\mathfrak{F}_\mathfrak{g}(V)$. Then we can verify the continuity of~$\rho_V$ similarly to the $C^\infty$-case using Lemma~\ref{l6.4} instead of Lemma~\ref{l5.2} and prove the topological injectivity in the same way as in Theorem~\ref{t5.1} taking into account the fact that all auxiliary assertions in \S\,\ref{ss4.5} hold in the case of the field~$\mathbb{C}$.
\end{proof}

We now turn to sheaves and introduce a multiplication in $\mathfrak{F}_\mathfrak{g}(V)$ in the general case.
Recall that a \textit{sheaf on a base of a topology} is defined similarly to a sheaf on a topological space, with the difference that all open sets participating in the axioms are assumed to belong to the base.
Denote by $\mathcal{B}$ the set of open subsets $V$ of~$\mathbb{C}^k$ such that polynomials are dense in~$\mathcal{O}(V)$. Since the polydiscs form a prebase of the topology on $\mathbb{C}^k$ and each polydisc is in $\mathcal{B}$, we conclude that $\mathcal{B}$ is a base of this topology. Thus we can consider sheaves on~$\mathcal{B}$. Identifying the Gelfand spectrum of~$\mathfrak{F}_\mathfrak{g}$ with~$\mathbb{C}^k$ we can assume that $\mathcal{B}$ is a base of the topology on the Gelfand spectrum. In what follows, we denote by $\tau_{VW}$ the  restriction map $\mathfrak{F}_\mathfrak{g}(V)\to \mathfrak{F}_\mathfrak{g}(W)$ for open subsets $V$ and $W$ of~$\mathbb{C}^k$ such that $W\subset V$.

\begin{theorem}
\label{t6.5}
Let $\mathfrak{g}$ be a nilpotent complex Lie algebra. Then the correspondences
\begin{equation}
\label{eq6.1}
V\mapsto \mathfrak{F}_\mathfrak{g}(V)\quad\text{and}\quad (W\subset V)\mapsto \tau_{VW}
\end{equation}
determine a sheaf $\mathfrak{F}_\mathfrak{g}(-)$ of Fr\'echet--Arens--Michael algebras on the base~$\mathcal{B}$.
\end{theorem}

The proof of Theorem~\ref{t6.5} is the same as for Theorem~\ref{t5.5} with the obvious replacement of the sheaf of smooth functions by the sheaf of holomorphic functions.

\begin{theorem}
\label{t6.6}
Let $\mathfrak{g}$ be a nilpotent complex Lie algebras.  Then for every open subset $V$ of~$\mathbb{C}^k$ there is a multiplication making  $\mathfrak{F}_\mathfrak{g}(V)$ a Fr\'echet--Arens--Michael algebra. Moreover, $\mathfrak{F}_\mathfrak{g}(V)$ is independent of the choice of a basis of~$\mathfrak{n}$. In the case when polynomials are dense in~$\mathcal{O}(V)$, the multiplication can be obtained by extension from~$U(\mathfrak{g})$ by continuity. Besides, the correspondences in~\eqref{eq6.1} determine a sheaf $\mathfrak{F}_\mathfrak{g}(-)$ of Fr\'echet--Arens--Michael algebras on the Gelfand spectrum of~$\mathfrak{F}_\mathfrak{g}$.
\end{theorem}

\begin{proof}
By Theorem~\ref{t6.5}, we have a sheaf of Arens--Michael algebras
on the base~$\mathcal{B}$. Using standard methods of sheaf theory (see, for example, \cite[Part~1, \S\,6.30]{39}) we can show that it  extends uniquely to a sheaf~$\mathfrak{F}'_\mathfrak{g}(-)$ of Arens--Michael algebras on~$\mathbb{C}^k$. On the other hand, arguing in the same way as in the proof of Theorem~\ref{t5.5}, we get that the correspondences in~\eqref{eq6.1} determine a sheaf $\mathfrak{F}_\mathfrak{g}(-)$ of Fr\'echet spaces on~$\mathbb{C}^k$, which is evidently is an extension of the sheaf on~$\mathcal{B}$. Note that $\mathfrak{F}'_\mathfrak{g}(-)$ and $\mathfrak{F}_\mathfrak{g}(-)$ are sheaves of locally convex spaces. Then, since the extension of a sheaf on a base of the topology to a sheaf on the whole space is unique,  we have  an isomorphism  $\mathfrak{F}'_\mathfrak{g}(V)\cong\mathfrak{F}_\mathfrak{g}(V)$ of locally convex spaces for every open subset~$V$. So we obtain the required multiplication in~$\mathfrak{F}_\mathfrak{g}(V)$. The rest is clear.
\end{proof}

\begin{remark}
\label{r6.7}
The explicit form of the multiplication in~$\mathfrak{F}_\mathfrak{g}(V)$ for an arbitrary open subset $V$ in~$\mathbb{C}^k$ is as follows. Let $(V_i)$ be a cover of~$V$ by elements of~$\mathcal{B}$. Then the left arrow in the diagram
$$
\mathfrak{F}_\mathfrak{g}(V) \to \prod_i \mathfrak{F}_\mathfrak{g}(V_i)\rightrightarrows \prod_{i,j} \mathfrak{F}_\mathfrak{g}(V_i\cap V_j)
$$
is an equalizer in the category of Arens-Michael algebras. So we can identify  $\mathfrak{F}_\mathfrak{g}(V)$ with a closed subalgebra  of a product of Arens-Michael algebras, the multiplication in each of which extends from~$U(\mathfrak{g})$ by continuity.
\end{remark}

\begin{remark}
\label{r6.8}
Comparing Theorem~\ref{t6.6} with Theorems~\ref{t5.1} and~\ref{t5.5}, we see that the statements of the last two are stronger since they provide not only Fr\'echet--Arens--Michael algebras but algebras with an additional property,  polynomial growth. While there is no such improvement in the theorems in this section. However, the shortcoming can be easily eliminated. Just as all $C^\infty_\mathfrak{g}(V)$ (in particular, $C^\infty_\mathfrak{g}$) are projective limits of Banach algebras of polynomial growth, so all $\mathfrak{F }_\mathfrak{g}(V)$ (in particular, $\mathfrak{F}_\mathfrak{g}$) are projective limits of Banach PI-algebras (those that satisfy a polynomial identity). Indeed, the property of being a PI-algebra is preserved under passing to a subalgebra and it is easy to see that $\mathcal{O}(V,\mathrm{T}_{p})$ and as well as its concomitant Banach algebras satisfy all the polynomial identities that $\mathrm{T}_{p}$ satisfies. The fact that $\mathrm{T}_{p}$ satisfies a polynomial identity is well known (see Remark~\ref{r2.10}).

Note also that every $\mathbb{R}$-algebra of polynomial growth is a projective limit of Banach PI-algebras (see Remark~\ref{r2.10} again).
\end{remark}

In conclusion, it is worth noting  that what is said in~\S\,\ref{ss5.2} about sheaves of non-commutative $C^\infty$-functions in the case when $\mathfrak{g}$ is not nilpotent is also true for non-commutative holomorphic functions.

\medskip

The author would like to thank D.~Belti\c{t}\u{a}, A.\,V.~Domrin, A.~Dosi, A.\,Yu.~Pirkovskii and Yu.\,V.~Turovski for useful consultations. The author is also grateful to the referees for valuable comments.

\end{document}